\documentclass{amsart}
\usepackage[left=4cm,right=4cm,top=3cm,bottom=3cm]{geometry}
\usepackage{amsfonts, amsthm, amssymb, amsmath, stmaryrd}
\usepackage{mathrsfs,array}
\usepackage{eucal,color,times,enumerate,accents}
\usepackage{tikz-cd}
\usepackage{bbm}
\usepackage[utf8]{inputenc}
\usepackage{pgfplots}
\pgfplotsset{compat=newest}
\usepackage{lipsum}
\usepackage{graphicx}
\usepackage{rotating}
\usepackage{enumitem}
\usepackage{enumerate}
\usepackage{footmisc}
\usepackage{textcomp}
\usepackage[bookmarksnumbered,bookmarksopen]{hyperref}
\hypersetup{colorlinks, linkcolor=black, citecolor=black}
\usepackage{cleveref}
\usepackage[colorinlistoftodos]{todonotes}
\usepackage[left=4cm,right=4cm,top=3cm,bottom=3cm]{geometry}
\usepackage{mathtools}
\newtheorem{theorem}[equation]{Theorem}

\newtheorem{lemma}[equation]{Lemma}

\newtheorem{corollary}[equation]{Corollary}

\numberwithin{equation}{subsection} 

\theoremstyle{definition}
\newtheorem{definition}[equation]{Definition}
\newtheorem{question}[equation]{Question}

\newtheorem{remark}[equation]{Remark}
\newtheorem{example}[equation]{Example}

\usepackage{tikz}
\usepackage{tikz-cd}
\usepackage{rotating}
\usetikzlibrary{arrows,shapes,positioning}
\usetikzlibrary{decorations.markings}
\tikzstyle arrowstyle=[scale=1]
\usepackage{lipsum}

\usepackage{url}
\usepackage{imakeidx}
\makeindex[title=Index of definitions and notations,intoc]
\newcommand{\Z}{\mathbb{Z}}
\newcommand{\Q}{\mathbb{Q}}

\newcommand{\rank}{\textrm{rank}}
\newcommand{\Dim}{\textrm{dim}}

\newcommand{\Spec}{\mathrm{Spec}}

\newcommand{\Hom}{\mathrm{Hom}}

\newcommand{\Gr}{\mathrm{Gr}}

\newcommand{\R}{\mathbb{R}}
\newcommand{\C}{\mathbb{C}}
\begin{document}
\title{Betti numbers of real semistable degenerations via real logarithmic geometry}
%\date{\today}
\author{Emiliano Ambrosi and Matilde Manzaroli}
\begin{abstract}
Let $X\rightarrow C$ be a totally real semistable degeneration over a smooth real curve $C$ with degenerate fiber $X_0$.
Assuming that the irreducible components of $X_0$ are simple from a cohomological point of view, we give a bound for the individual Betti numbers of a real smooth fiber near $0$ in terms of the complex geometry of the degeneration. This generalizes previous work of Renaudineau-Shaw, obtained via combinatorial techniques, for tropical degenerations of hypersurfaces in smooth toric varieties. The main new ingredient is the use of real logarithmic geometry, which allows to work with not necessarily toric degenerations.
\end{abstract}
\maketitle
\tableofcontents
\section{Introduction}
Let $X$ be a real algebraic variety, let $X(\mathbb C)$ be the set of its complex points and $X(\mathbb R)$ the set of its real points. For a topological space $Y$, set $b_i(Y):=\Dim_{\mathbb F_2}(H^i(Y,\mathbb Z/2\Z))$ for its $i^{th}$ Betti number.
\subsection{Individual Betti numbers}
The Smith-Thom inequality 
\begin{equation}
\label{eqn: smiththominequality}
	\sum_{i}b_{i}(X(\mathbb R))\leq \sum_{i}b_i(X(\mathbb C))
\end{equation}
bounds the total Betti number of the real points of $X$ with the total Betti number of its complexification. A variety is said to be \textit{maximal}, in the sense of Smith-Thom, if (\ref{eqn: smiththominequality}) is an equality.
The problem of finding (non trivial) bounds for the individual Betti numbers  $b_i(X(\mathbb R))$ in terms of the geometry of $X(\mathbb C)$ is a central topic in real algebraic geometry. For example, it was conjectured by Viro in \cite{Viro80} that if $X$ is a smooth projective real surface such that $X(\mathbb C)$ is simply connected then $b_1(X(\mathbb R))\leq \Dim(H^{1}(X,\Omega^1_X))$. Even if the conjecture has been disproved by Itenberg in \cite{Ite93, Iten97}, many smooth projective real algebraic varieties constructed with the tools at our disposal (for example, using Viro's patchworking \cite{Viro83,Viro84}) verify the inequality
\begin{equation}\label{subhodgemaximal}
	b_i(X(\mathbb R))\leq\sum_{j} h^{i,j}(X)
\end{equation}
where $h^{i,j}(X):=\dim(H^i(X,\Omega^j_X))$. Following \cite{BS22}, if a smooth projective real variety has no torsion in cohomology and satisfies the inequality (\ref{subhodgemaximal}), we say that it is \textit{sub-Hodge expressive}.

The main goal of this paper is to explain this phenomenon in some cases. Roughly speaking, the general principle is that, often, the constructions of real varieties with prescribed topology are achieved first by constructing a degenerated version of the pair $(X(\mathbb C),X(\mathbb R))$, whose irreducible components are simpler to deal with and then by deforming them back in non-trivial ways. Our main results (Theorem \ref{main}, Corollary \ref{cor: main}), shows that if these pieces are very simple, from a cohomological point of view, then the variety obtained gluing these pieces satisfies the inequality (\ref{subhodgemaximal}), up to some torsion. The main novelty of our approach is to use (real) logarithmic geometry, which allows to extend previous results for toric degenerations to more general, non necessarily toric, families.
\subsection{The case of tropically smooth hypersurfaces}\label{tropicalintro}
\subsubsection{Renaudineau-Shaw theorem}\label{subsubsec: Artie_Shhaw_result}
The most general result previously known in this setting is \cite[Theorem 1.4]{RenShaw22}, which proves a conjecture of Itenberg (\cite{Iten17}) and whose proof involves real and complex combinatorial algebraic geometry arguments. Renaudineau and Shaw show that for a real compact hypersurface $X$ near the $\Q$-regular smooth tropical limit inside a smooth toric variety (see \cite{tropicalhomology} for the definitions involved), one has the following inequality (\cite[Theorem 1.5]{RenShaw22}):
$$b_i(X_t(\mathbb R))\leq \sum_{j}\dim(H_{i,j}(Trop(X),\mathbb Z/2\Z)),$$
where $H_{i,j}(Trop(X),\mathbb Z/2\Z)$ is the $(i,j)$ tropical homology group with $\mathbb Z/2\Z$-coefficients (which is independent from the real structure of the family) and $X_t(\mathbb R)$ are the real points of a smooth fiber near the tropical limit. Moreover, by \cite{ArnRenSha21}, one knows that the variant $H_{i,j}(Trop(X),\mathbb Z)$ of tropical homology with integer coefficients is torsion free, while by \cite{tropicalhomology} one knows that its rank is the $(i,j)$ Hodge number of the general fiber $X_t$ of the family. Summing up, one gets the sub-Hodge expressive inequality (\cite[Theorem 1.4]{RenShaw22}):
$$b_i(X_t(\mathbb R))\leq \sum_{j}\dim(H^{i}(X_t,\Omega^j)).$$

In this paper we generalize this result to more general, non necessarily  toric, families of varieties, avoiding the use of combinatorial arguments and constructing, via real logarithmic geometry, a space which allows to simultaneously relate real and complex data.
\numberwithin{equation}{subsubsection} 
\begin{remark}
Very recently, there have been other two generalizations of the results and the techniques in \cite{RenShaw22}:
	\begin{itemize} 
		\item On the combinatorial side, in \cite{BLMR22}, Brugall\'e, L\'opez de Medrano and Rau extend the main results of \cite{RenShaw22} to $PL$-manifolds associated to non-necessarily convex triangulations of non-singular lattice polytopes. 
		\item On the geometric side, in a work in progress of Renaudineau, Shaw and Rau, the authors generalize \cite[Theorem 1.5]{RenShaw22} and the strategy therein to general families with smooth tropicalisation near the tropical limit inside a smooth toric variety.
	\end{itemize}
	\end{remark}
\subsubsection{Geometry of tropicalisation}\label{subsec : geometrytrop}
In order to put our main results (Theorem \ref{main}, Corollary \ref{cor: main}) and their hypothesis into context, we present an overview of the algebraic geometry behind the tropicalisation. In the setting of Section \ref{subsubsec: Artie_Shhaw_result}, one can construct from $Trop(X)$ a toric degeneration $T \rightarrow C$ of the ambient toric variety and an associated one dimensional toric degeneration $X\rightarrow C$ of $X$, whose intersection with every toric orbit of $T$ is the complement of an hyperplane arrangement endowed with the standard complex conjugation. Complements of hyperplane arrangements endowed with the standard complex conjugation have tight properties, both from the complex and the real point of view:
\begin{enumerate}[label=(\alph*)]
\item the only non trivial cohomology group of their real points is the $0$-th;
\item they are maximal varieties, in the sense of (\ref{eqn: smiththominequality});
\item the mixed Hodge structure on the $i^{th}$-cohomology group with rational coefficients is pure of type $(i,i)$.
\end{enumerate}
In the proof of the results of \cite{tropicalhomology}, property (c) is used explicitly to give a geometric interpretation to the combinatorial complex computing $H_{i,j}(Trop(X),\mathbb Q)$ (see \cite[Section 5.3]{tropicalhomology}), while in \cite{RenShaw22} properties (a) and (b) are used implicitly to compare the dimensions of different combinatorial vector spaces associated to $Trop(X)$; see \cite[Corollary 3.15]{RenShaw22}.
\subsection{Main result}
\numberwithin{equation}{subsection} 
Assume that $C$ is a smooth real curve and $f:X\rightarrow C$ a real projective morphism  which is smooth outside a real point $0\in C(\mathbb R)$ and strictly-semistable around $0$, in the sense that the irreducible components of $X_0$ are smooth and, locally analytically around $0$, the family $f:X(\mathbb C)\rightarrow C(\mathbb C)$ is isomorphic to the standard semistable degeneration $\Spec(\mathbb C[x_1,\dots, x_n,T]/(x_1\dots x_n-T))\rightarrow \Spec(\mathbb C[T]).$ Assume furthermore that $f:X\rightarrow C$ is totally real, i.e. that the irreducible components of $X_0(\mathbb C)$ are real.
Write 
$$X_0=\bigcup_{i\in I}X_i$$
 for the decomposition of $X_0$ in irreducible components and for every subset $J\subseteq I$ set 
 $$X_J:=\bigcap_{i\in J}X_i\quad \text{and} \quad X^{0}_J:=X_J\setminus \bigcup_{i\not \in  J}X_i.$$
Then 
 $$X_0=\coprod_{J\subseteq I}X^0_J$$
 is a stratification $\mathfrak I:=\{X^0_J\}_{J\subseteq I}$ of $X_0$ by smooth real algebraic subvarieties. Fix a refinement $\mathfrak Z:=\{X^0_\Delta\}$ of $\mathfrak I$, made of smooth real algebraic varieties. 
 
 In Section \ref{subsec: construction} we construct, for every ring $A$ and every $q\geq 0$ a canonical cochain complex $C_{q,\mathfrak Z,A}^{\bullet}$ of $A$-modules depending only on the complex geometry of the stratification $\mathfrak Z$. This complex plays, in this general setting, the role of the combinatorial complex computing $H_{i,j}(Trop(X),A)$ used in \cite{RenShaw22, ArnRenSha21, tropicalhomology}. 
 Inspired by the discussion in Section \ref{subsec : geometrytrop}, we consider the following condition on the members of  $\mathfrak Z$.
\begin{enumerate}
	\item[(a)] $H^i(X^0_\Delta(\mathbb R),\mathbb Z/2\Z)=0$ for all $i\geq 1$ and $X^0_\Delta\in \mathfrak Z$;
	\item[(b)]  $X^0_\Delta$ is a maximal variety, for all $X^0_\Delta\in \mathfrak Z$;
	\item[(c)] the mixed Hodge structure on $H^i(X^0_\Delta,\Q)$ is pure of type $(i,i)$ and $H^i(X^0_\Delta,\Z)$ is torsion free, for all $i\geq 1$ and $X^0_\Delta\in \mathfrak Z$.
\end{enumerate}
Our main result is the following.
\begin{theorem}\label{main}
\begin{enumerate}
\item[]	
	\item Assume that (a) and (b) hold. 
	Then, for every $t\in C(\mathbb R)$ close to $0$ one has:

 	$$b_p(X_t(\mathbb R))\leq \sum_{q}\Dim(H^p(C_{q,\mathfrak Z,\Z/2\Z}^{\bullet})).$$
	\item Assume that (a),(b) and (c) hold. Then for every $t\in C(\mathbb R)$ close to $0$, one has:
	\begin{enumerate}
	\item[(i)] $\Dim(H^p(C_{q,\mathfrak Z,\Z}^{\bullet}\otimes \Q))=h^{p,q}(X_t)$
	\item[(ii)] $C_{q,\mathfrak Z,\Z}^{\bullet}\otimes \Z/2\Z \simeq C_{q,\mathfrak Z,\mathbb Z/2\Z}^{\bullet}$
	\end{enumerate}
\end{enumerate}

 \end{theorem}

%
%\begin{theorem}\label{main}
%\begin{enumerate}
%\item[]	
%	\item Assume that 
% 	\begin{enumerate}
%\item[(a)] $H^i(X^0_\Delta(\mathbb R),\mathbb Z/2\Z)=0$ for all $i\geq 1$ and $X^0_\Delta\in \mathfrak Z$;
% 		\item[(b)]  $X^0_\Delta$ is a maximal variety, for all $X^0_\Delta\in \mathfrak Z$.
% 	\end{enumerate}
% 	Then, for every $t\in C(\mathbb R)$ close to $0$ one has:
% 	$$b_p(X_t(\mathbb R))\leq \sum_{q}\Dim(H^p(C_{q,\mathfrak Z,\Z/2\Z}^{\bullet})).$$
%\item
% 	In addition to (a) and (b), assume that
% 	\begin{enumerate}
% 		\item[(c)] 
% 	the mixed Hodge structure on $H^i(X^0_\Delta,\Q)$ is pure of type $(i,i)$ and $H^i(X^0_\Delta,\Z)$ is torsion free. 
% 	\end{enumerate} 
%Then, for every $t\in C(\mathbb R)$ close to $0$, one has:
%
%	$$(i)\quad \Dim(H^p(C_{q,\mathfrak Z,\Z}^{\bullet}\otimes \Q))=h^{p,q}(X_t)\quad \text{and} \quad (ii) \quad C_{q,\mathfrak Z,\Z}^{\bullet}\otimes \Z/2\Z \simeq C_{q,\mathfrak Z,\mathbb Z/2\Z}^{\bullet}.$$
%
%
% 	\end{enumerate}
% \end{theorem}
 Theorem \ref{main} directly implies the following corollary, which was the main motivation for this paper. 
 \begin{corollary}
 \label{cor: main}
 Assume that (a), (b), (c) hold and that $H^p(C_{q,\mathfrak Z,A}^{\bullet})$ is torsion free for every $p,q\in \mathbb N$. Then for every $t\in C(\mathbb R)$ close to $0$ and every $p\in \mathbb N$ one has
 $$b_p(X_t(\mathbb R))\leq \sum_{q}h^{p,q}(X_t).$$
 \end{corollary}
 \begin{remark}
 	We point out that even in very easy examples is necessary to refine the standard stratification $\{X^0_J\}$ in order to guarantee that the strata satisfy the hypothesis in Theorem \ref{main}. If ones take the trivial degeneration of $\mathbb P^1$, the stratification $\{X^0_J\}$ is the trivial one and hence it does not satisfy the hypothesis. So, one has to further stratify $\mathbb P^1$, for example as
 	 	$$\mathbb P^1=\Big(\mathbb P^1-\Big (\{[0,1]\}\coprod \{[1,0]\}\Big)\Big)\coprod \{[0,1]\}\coprod \{[1,0]\},$$
 	which corresponds to the stratification of $\mathbb P^1$ into toric orbits.
 \end{remark}
In the case of degenerations constructed from a smooth $\mathbb Q$-regular tropicalisation inside a smooth toric variety, the stratification $\mathfrak Z:=\{X^0_\Delta\}$ is the stratification obtained intersecting the special fiber with the toric orbits of the ambient toric variety, hence we see that our results generalise those in \cite{RenShaw22} to non necessarily toric degenerations. 

Under the hypothesis of Theorem \ref{main}, in order to understand the link between the Hodge numbers and the real Betti numbers, it remains to understand the torsion in the cohomology of $C_{p,\mathfrak Z,\mathbb Z}^{\bullet}$. While torsion appears in general, this raises the following question.
\begin{question}
	Let $X$ be a smooth complete intersection having smooth $\Q$-regular tropicalisation inside a smooth toric variety. Is the cohomology of $C_{p,\mathfrak Z,\mathbb Z}^{\bullet}$ torsion free?
\end{question} 
\subsection{Strategy}\label{sec: introstrategy}
\subsubsection{Main ideas}
In order to attack this kind of problems, the general involved strategy was to separately compute complex and real information and, in a second moment, show that they are related. For example, in \cite{tropicalhomology} it is shown that tropical homology is related to Hodge numbers, in \cite{RenShaw22} it is first shown that the tropical interpretation of the Viro patchworking method (\cite{Viro83, Viro84}) gives a way to compute the real Betti numbers via the (co)homology of a tropical sheaf on the tropical varieties (\cite[Theorem 3.7  and Remark 3.8]{RenShaw22}) and then it is shown how to relate the (co)homology of this sheaf to tropical homology (\cite[Section 4]{RenShaw22}).

 A different approach to similar problems has been recently proposed by Brugall\'{e} in \cite{Bru22}. There, complex and real invariants are shown to satisfies the same gluing relations under totally real semistable degenerations, so that, in order to relate the global invariants, it is enough to relate the local invariants, which are easier to compute. Inspired by this, our basic strategy is to relate the real Betti and the Hodge numbers via the geometry of a common ambient space. The main innovation of this paper is the use of real logarithmic geometry to construct and study this common ambient space, which allows to use a more sophisticated and less combinatorial machinery. After the construction of such a space, realised in Sections \ref{subsec: semistablecasean} and \ref{subsec : semistablecase real}, the cohomology of the general real fiber can be computed by a filtered complex; see Section \ref{subsec : constructionfiltration}. The idea of the use of filtered complexes is inspired by \cite{RenShaw22}, where it was constructed via combinatorial techniques. Since these combinatorial approach are not available in our general setting, we use a different approach based on equivariant cohomology.
  \begin{remark}
 	Observe that \cite[Proposition 2.1]{Bru22} applies under the hypothesis of Theorem \ref{main} and implies that, for all $t \in C(\mathbb R)$ close to $0$, one has $\chi= \sigma$, where $\chi$ is the Euler characteristic of $X_t(\mathbb R)$ and $\sigma$ is the signature of $X_t(\mathbb C)$.
 \end{remark}
 \subsubsection{Equivariant cohomology}
 \label{introequivariant}
 To explain our strategy in more details, let us recall a modern proof of the Smith-Thom inequality (\ref{eqn: smiththominequality}), as in \cite[Chapter 4, IV, Pag. 55]{borel}. Let $G$ be $\mathbb Z/2\Z$ acting via the complex conjugation on $X$. Then there is a spectral sequence
 $$E_2^{a,b}:=H^a(G,H^b(X,\mathbb Z/2\Z))\Rightarrow H^{a+b}_G(X,\mathbb Z/2\mathbb Z)$$
 where $H^a(G,H^b(X,\mathbb Z/2\Z))$ is the group cohomology of $G$ acting on $H^b(X,\mathbb Z/2\Z)$ and $H^{a+b}_G(X,\mathbb Z/2\mathbb Z)$ is the $G$-equivariant cohomology of $X$. Since  $H^{a+b}_G(X,\mathbb Z/2\mathbb Z)\simeq H^*(X(\mathbb R),\mathbb Z/2\mathbb Z)$ for $a+b>2\Dim(X)$, for every $m>2\Dim(X)$ one gets a filtration $_mF^i$ of $H^*(X(\mathbb R),\mathbb Z/2\mathbb Z)$ such that $_mF^i/ _mF^{i+1}$ is a subquotient of $H^i(G,H^{m-i}(X,\mathbb Z/2\Z))$. Since $\Dim(H^a(G,H^b(X,\mathbb Z/2\Z)))\leq \Dim(H^b(X,\mathbb Z/2\Z))$, one gets the desired inequality. This argument shows also that $X$ is maximal if and only if the spectral sequence degenerates at $E_2$ and the action of $G$ on $H^*(X,\mathbb Z/2\mathbb Z)$ is trivial. Since we assume maximality for every stratum of our degeneration, we could try to apply this argument to each stratum and then glue them together to relate complex and real geometry of the general fiber; see Section \ref{subsec : equivariant} for more details. 
\subsubsection{Real logarithmic geometry}
From a topological point of view, as remarked in \cite{Bru22}, the general real fiber is the union of coverings of the real strata of the special fiber. Naively, one can hope to construct a stratification of the real general fiber from the one of the real special fiber gluing this covering.
Unfortunately, it is unclear to us how to compare this kind of constructions with the complex geometry of the general fiber, since the gluing conditions might be quite complicated. On the other hand, what one can do is to use real and complex logarithmic geometry to simultaneously stratify the complex and the real general fiber; see Sections \ref{subsec : loog} and \ref{sec: real_log}.

More precisely, there is a natural structure  $X^{log}$ and $C^{log}$ of real-log variety, in the sense of \cite{Argu21},  on $X$ and $C$ making the morphism  $X^{log}\rightarrow C^{log}$ a smooth morphism of real log-varieties. Taking the fiber in $0$, we get a morphism $X_0^{log}\rightarrow 0^{log}$. Then we can take the analytification, in the sense of Kato-Nakayama \cite{katonakayama}, of this morphism to construct a morphism of $C^{\infty}$-manifolds with corners $X_0^{log}(\mathbb C)\rightarrow S^1$ endowed with involutions. The main point of logarithmic geometry is that the fiber $X_{0,1}^{log}(\mathbb C)$ at $1$ of this morphism is homeomorphic to the general fiber of the family (\cite{nakayamaogus}), while its real part $X_{0,1}^{log}(\mathbb R)$ is homeomorphic  to the real part of a general (positive) fiber (\cite{Argu21} or \cite{Rau22}). This $X_{0,1}^{log}(\mathbb C)$ is the common ambient space we were looking for. 
\begin{figure}[h!]
\begin{center}
\begin{picture}(100,250)
\put(-125,-150){\includegraphics[width=0.8\textwidth]{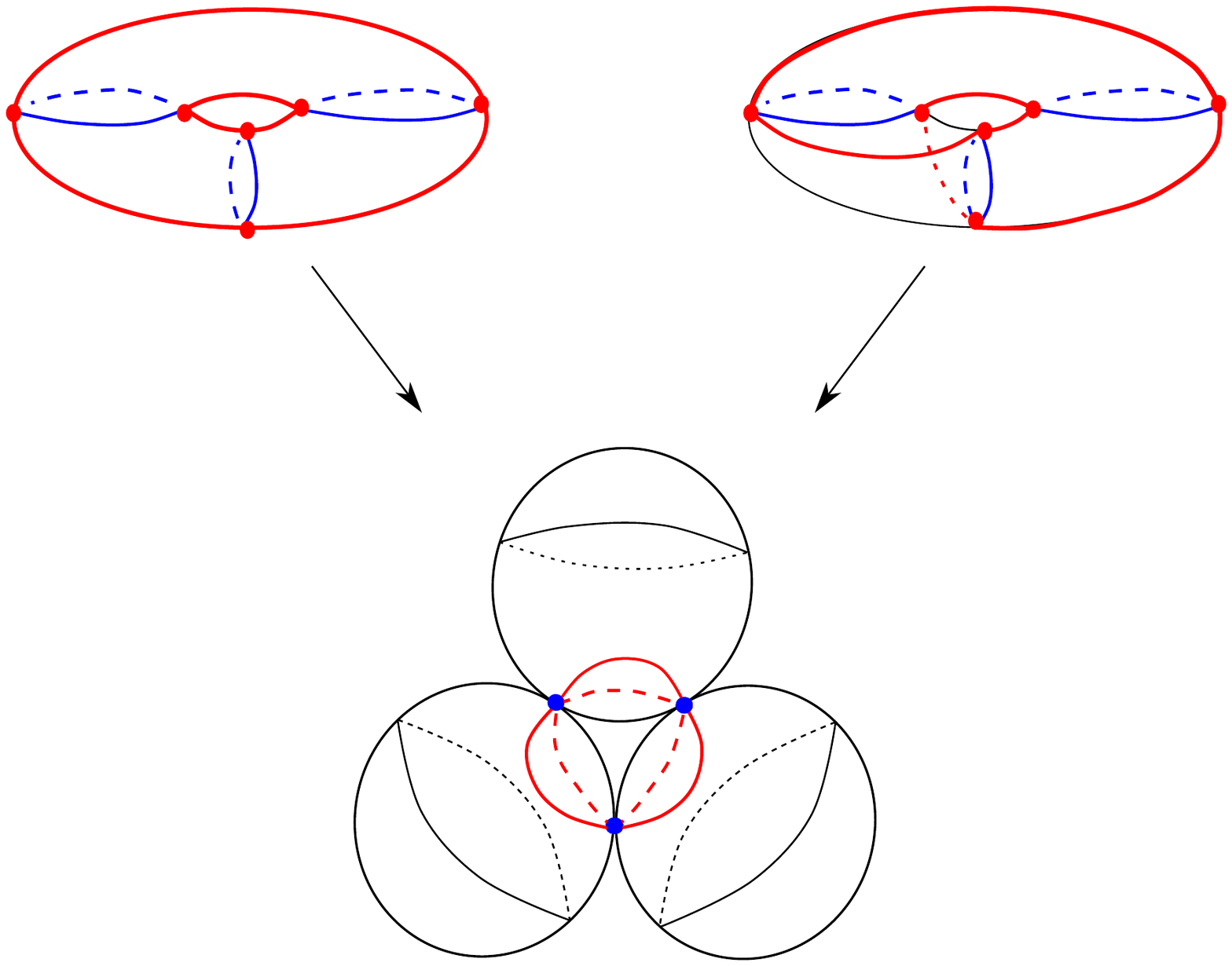}}
\end{picture}
\end{center}
\caption{ $(X^{log}_{0,1}(\mathbb C),X^{log}_{0,1}(\R)) \rightarrow (X_0(\mathbb C),X_0(\mathbb R))$. On the left a maximal case and on the right a non-maximal one.}
\label{fig: real_logan_cubic}
\end{figure}
%dopo disegni
\textcolor{white}{a}

It remains to relate this space to the special fiber. The formalism of logarithmic geometry gives a morphism $X_{0,1}^{log}(\mathbb C)\rightarrow X_{0}(\mathbb C)$ compatible with the involutions, which is a nice locally constructible fibration with fiber at $x$ equal to $S^1$ at some power $n_x$, depending only on which stratum $X^0_J$ the point $x$ lies. This allows to relate the complex geometry of the general fiber to that of the special fiber. To understand the real part, one has just to take the fixed points of the involutions, and thanks to the totally real assumption, one gets a locally constructible cover $\pi: X_{0,1}^{log}(\mathbb R)\rightarrow X_{0}(\mathbb R)$ with fiber at $x$ equal to $\{\pm 1\}$ at some power $n_x$ depending only on which stratum $X^0_J$ the point $x$ lies; see Section \ref{subsec stratificationacaso}. The situation is explained in an example in Fig.\ref{fig: real_logan_cubic}, in the case of a family of real elliptic curves (endowed with two different real structures, respectively on the left and on the right of Fig.\ref{fig: real_logan_cubic}) degenerating to the union of three genus $0$ curves. 

\subsubsection{Filtred complexes}
In \cite{RenShaw22}, in order to prove the main theorem, it is given an explicit complex which computes combinatorially the real Betti numbers and which is filtered by subcomplexes $K_{i}$, in way that the resulting graded quotients are isomorphic to the complex computing the tropical homology. The bound is then obtained as a direct consequence of the spectral sequence for a filtered complex.
Taking inspiration from the above argument, we now have all the ingredients to prove the desired bound.
Once $X^{log}_{0,1}(\mathbb C)$ and $X^{log}_{0,1}(\mathbb R)$ are stratified in a compatible way via the stratification of the special fiber, one can run the argument of Section \ref{introequivariant}, using equivariant cohomology on each stratum. First, since the real part of the strata have only $H^0$, the spectral sequence computing the cohomology of a stratified space reduces to a single complex, so that we can compute the cohomology groups of $X^{log}_{0,1}(\mathbb R)$ via a complex (Section \ref{computationrealnumber}) whose terms are the cohomology of the strata. Then, since the strata are maximal and have only $H^0$, the spectral sequence for equivariant cohomology gives a filtration of the real cohomology of each real stratum whose graduated complex is the cohomology of the corresponding complex stratum. Putting these together, we get the filtration that we were looking for and we conclude by using the spectral sequence of a filtered complex. This ends the proof of Theorem \ref{main}(1).
\subsubsection{Hodge numbers}
We finally explain how to get Theorem \ref{main}(2). By the arguments in \cite[Pag. 31]{tropicalhomology}, the theory of limiting Hodge structures and Theorem \ref{main}'s hypothesis, one has that the $(p,q)$-Hodge numbers of the general fiber can be computed as dimension of the weight $2q$ part in the limiting Hodge structure of $X^{log}_{0,1}(\mathbb C)$ of degree $p+q$. Hence it is enough to show that  the cohomology of the complex $C_{q,\mathfrak Z,\Z}^{\bullet}\otimes \Q$ computes the weight filtration on the limiting Hodge structure. 

While this would follows easily from a Hodge theory for the cohomology with compact support of open logarithmic varieties, this theory seems not to be developed. One of the problem is that the open stratum $X^0_{J}$ endowed with its natural log-structure is not a log-smooth variety, so that one can not naively appeal to some form of logarithmic Poincar\'e duality.

To avoid this problem, we filter $C_{q,\mathfrak Z,\Z}^{\bullet}\otimes \Q$ by exploiting the Leray spectral sequence for the morphisms $X^{log}_{\Delta,1}(\mathbb C)\rightarrow X^0_{\Delta}(\mathbb C)$, to relate the cohomology of $C_{q,\mathfrak Z,\Z}^{\bullet}\otimes \Q$ to the one of the nearby cycles sheaves on $X^0_{\Delta}(\mathbb C)$. At the same time, we compute the weight filtration on $H^n(X^{log}_{0,1}, \mathbb Q)$ by exploiting the Leray Spectral sequence for the morphism $X^{log}_{0,1}(\mathbb C)\rightarrow X_0(\mathbb C)$, also called the vanishing cycles spectral sequence. The main point is that the assumption on the weights allows us to show that the terms appearing in the different spectral sequences involved coincide, since they both compute the same weight part of a mixed Hodge structure. Since the details are somehow involved, we refer the reader to Section \ref{sec : strategy}, for a more precise account.
\subsection{Organisation of the paper}

The paper is organised as follows. In Section \ref{notation}, we fix some basic notation. Section \ref{subsec : loog} and \ref{sec: real_log} are devoted to recall and complement some properties of complex and real log-varieties respectively. In Section \ref{sec: equivariant_coho}, we prove the first point of Theorem \ref{main} and in Section \ref{sec: hodge_number} we prove the second point of Theorem \ref{main}.
	\subsection{Acknowledgement}

	The authors are grateful to Erwan Brugall\'e, Taro Fujisawa, Kris Shaw, Atsushi Shiho for helpful discussions. A thanks to Michele Ancona, Erwan Brugall\'e and Olivier Benoist for useful comments on a preliminary version of this paper. The authors thank Athene Grant and CAS young fellows program "Real structures in discrete, symplectic, and tropical geometries". 
\section{Notation and conventions}
\numberwithin{equation}{subsubsection} 

\subsection{Stratification}	\label{notation}
\numberwithin{equation}{subsubsection} 

Let $X$ be a topological space and $\{X^0_{\Delta}\}$ a collection of locally closed subspaces of $X$. We write  $\{X_{\Delta}\}$ for the closure of $X^0_{\Delta}$. We say that $\{X^0_{\Delta}\}$ is a stratification of $X$ if $$X=\coprod X^0_{\Delta} \quad \text{and for every $\Delta$}  \quad X_{\Delta}=\coprod_{X^0_{\Delta'}\subseteq X_{\Delta}}X^0_{\Delta'}.$$

Let $\{X^0_{\Delta'}\}$ another stratification of $X$. We say that $\{X^0_{\Delta'}\}$ is a refinement of $\{X^0_{\Delta}\}$ if every $X^0_{\Delta'}$ is included is some $X^0_{\Delta}$ and for every $X^0_{\Delta}$ the collection $\{X^0_{\Delta'}\}_{X^0_{\Delta'}\subseteq X^0_{\Delta}}$ is a stratification of $X^0_{\Delta}$.

\subsection{Degenerations}
Let $C$ be a curve over a field $k$, a point $0\in C(k)$ and $f:X\rightarrow C$ a projective morphism of smooth varieties. We say that $f$ is semistable (around $0$) if it is smooth over $C-\{0\}$ and if, \'{e}tale locally around $0$, it is isomorphic to the standard semistable degeneration $\Spec(k[x_1,\dots,x_n,T]/(x_1\dots x_n))\rightarrow \Spec(k[T])$. We say that $f:X\rightarrow C$ is strictly semistable (around $0$) and the irreducible components $\{X_i\}_{i \in I}$ of the fiber $X_0$ in $0$ are smooth. Setting for every subset $J\subseteq I$
 $$X_J:=\bigcap_{i\in J}X_i\quad \text{and} \quad X^{0}_J:=X_J\setminus \bigcup_{i\not \in  J}X_i,$$
the collection $\mathfrak I:=\{X^0_J\}_{J\subseteq I}$ is a stratification of $X_0$ by smooth algebraic subvarieties. 
If $k=\mathbb R$, we say that $f$ is totally real, if the irreducible components of $X_0(\mathbb C)$ are real.
\subsection{Local orientations}\label{sec : orientation}
Let $C$ be a smooth real curve and $0\in C(\mathbb R)$ a real point. The set $C(\R)$ is  homeomorphic to a disjoint union of open subsets of $S^1$'s. Denote with $\mathcal H$ the real connected component of $C(\R)$ containing $0$. We say that $t\in  C(\mathbb R)$ is \textit{near} $0$, if $t$ belongs to $\mathcal H$. Let $D\subseteq \mathbb C$ be an open disc centered in $(0,0)$ and endowed with the canonical involution induced by the complex conjugation acting on $\C$.
	    	
	 An \textit{orientation around} $0 \in C$ is a pair $(U,\psi)$ where $U$ is an open neighbourhood of $0$ in $C(\mathbb C)$ stable under complex conjugation and $\phi$ is an equivariant homeomorphism $\psi:U\rightarrow D$ sending $0$ to $0$. An orientation induces an homeomorphism $\phi_{\mathbb R}:U(\mathbb R)\xrightarrow{\simeq}(-1,1)$, where $U(\mathbb R)$ is the set of points of $U$ fixed by the involution. Let $U_{+}:=\phi_{\mathbb R}^{-1}((0,1))$. A point $t\in C( \mathbb R)$ near $0$ is said to be \textit{positive} with respect to an orientation $(U,\psi)$ if it lies in $(\mathcal H \cap U_+)\setminus\{0\}$; otherwise \textit{negative}. 
	 
	 Assume that $X\rightarrow C$ is a morphism of real algebraic varieties.
	  A fiber $f^{-1}(t)$, for $t\in C( \mathbb R)$, is said to be \textit{near} $f^{-1}(0)$ if $t$ is near $0$. Moreover, for some fixed orientation $(U,\phi)$ around $0$, a fiber $f^{-1}(t)$ is said to be \textit{positive} if $t$ is near $0$ and it is positive with respect to $(U,\psi)$; otherwise \textit{negative}. 

  \section{Log-varieties and analytification}\label{subsec : loog}
  In this section we construct the complex $C_{p,\mathfrak Z,A}^{\bullet}$ appearing in Theorem \ref{main} and we state Theorem \ref{main} in a more precise way (Theorem \ref{mainpreciso}). To this end, we need to recall a few basic notions from logarithmic geometry (Section \ref{subsec : loggeo}), and from their associated analytification  (Section \ref{subsec: analytification}). Then we study the geometry of certain stratifications of these analytifications (Section \ref{subsec: stratification}) and we finally construct the complex $C_{p,\mathfrak Z,A}^{\bullet}$ in Section \ref{subsec: construction}.
  \subsection{Logarithmic varieties}\label{subsec : loggeo}
  \subsubsection{Definitions and examples}

  We start recalling a few basic basic notions from logarithmic geometry. For more details see for example \cite{illusielog} and \cite{Kato89}.
  
  Let $k$ be a field and $X$ a $k$-variety. We consider $\mathcal O_X$ as a sheaf of (abelian) monoids with the multiplication operation.
  \begin{definition}
  	\begin{enumerate}
  	\item[]
  		\item A pre-log structure on $X$ is a morphism of sheaf of (abelian) monoids $\alpha_X:M\rightarrow \mathcal O_X$.
  		\item A pre-log structure $\alpha_X:M\rightarrow \mathcal O_X$ is said a log structure if the natural morphism $\alpha_{\vert \alpha^{-1}(\mathcal O_X^*)}:\alpha^{-1}(\mathcal O_X^*)\rightarrow \mathcal O_X^*$  is an isomorphism.   		
  		\item Given a pre-log structure  $\alpha_X:M\rightarrow \mathcal O_X$, there is a canonical associated log-structure $\alpha'_X:M\rightarrow \mathcal O_X$ (see \cite[Section 1.3]{Kato89}).
  		\item A log-variety is a pair $(X,M)$ where $X$ is a $k$-variety and $M$ a log-structure on $X$.
  	\end{enumerate} 
  \end{definition}
  All the logarithmic structures in this paper will be assumed to be fine and saturated (\cite[1.3]{illusielog}) and all the operation between log-varieties will be done in the category of fine and saturated log varieties. Given a log variety $(X,M)$, we consider $\mathcal O_X^*$ as a submonoid of $M$ via $\alpha$. We write $M^{gr}$ for the groupification of the sheaf of monoids $M$ and $\overline M$ (resp. $\overline M^{gr}$) for the quotient $M/\mathcal O_X^*$ (resp. $\overline{(M/\mathcal O_X^*)}^{gr})$. There is a natural notion of morphism between logarithmic varieties. Given a log-variety $(X,M)$ and a morphism of $k$-varieties $f:Y\rightarrow X$, there is an induced log-structure $f^*M$ on $Y$ (\cite[Section 1.4]{Kato89}) and a morphism of log varieties $f:(Y,f^*M)\rightarrow (X,M)$ (morphisms of this type are called strict). Here some examples that are fundamental for our purposes and, therefore, that will be treated further in Examples \ref{exampleanilytification}, \ref{examplerealanalityfication}.
  \begin{example}\label{examplelogstructure}
  In the following examples anytime a morphism of monoids $\mathbb N \rightarrow \mathbb N^r$ is given, it is the diagonal map which sends $1$ to $\sum \limits_{i=1}^r e_i$, where $e_i$ denotes the $n$-uple with the $i$-th component equals $1$ and the remaining $0$.
  	\begin{enumerate}
  	\item The trivial log-structure on $X$ is given by $(X,\mathcal O_X^*)$ and $\alpha_X$ is the natural inclusion $\mathcal O^*_X\rightarrow \mathcal O_X$. In this way any variety can be consider as a log-variety with the trivial log-structure. If $(X,M)$ is a log-variety, there is a canonical morphism $(X,M)\rightarrow X$.
  	\item Let $X=\Spec(k)$ and let $\mathbb N\rightarrow k$ be the map sending $1$ to $0$. The associated log-structure $\mathbb N \times k^{*} \rightarrow k$ sends $(1, \lambda)$ to $0$ and $(0,\lambda)$ to $\lambda$. It is called the standard log-structure on the point and we write $\Spec(k)^{\log}$ for the associated log-scheme. 
  	\item Let $C$ be $\mathbb A^1=Spec(k[T])$ and let $\mathbb N\rightarrow k[T]$ be the map sending $1$ to $T$. The associated log structure $\mathbb N\times k^*\rightarrow k[T]$ sends $(n,f)$ to $fT^n$. It is called the log structure on $\mathbb A^1$ associated to the divisor $0$. We write $C^{log}$ for the associated log-scheme. The restriction of such log-structure to $0\in \mathbb A^1=Spec(k[T])$, induces the standard log-structure on the point. 
   	\item Let $X=\Spec(k[x_1\dots x_n,T]/(x_1\dots x_n-T))$ and $\mathbb N^n\rightarrow k[x_1\dots x_n,T]/(x_1\dots x_n-T)$ be the pre-log structure sending $e_i$ to $x_i$. The associated log structure $\mathbb N^n\times k^*\rightarrow k[x_1\dots x_n,T]/(x_1\dots x_n-T)$ sends $(e_i,f)$ to $fx^i$. We write $X^{log}$ for the associated log-variety. If $C=\mathbb A^1$ is endowed with the log-structure in Example \ref{examplelogstructure}(3), i.e. the previous example, the natural morphism  $X\rightarrow C$, extends to a morphism of log varieties $X^{log}\rightarrow C^{log}$ described by the commutative diagram
   	   \begin{center}
			 \begin{tikzcd}
		\mathbb N \arrow{r}{1 \mapsto T} \arrow{d}{1 \mapsto \sum e_i}& k[T]\arrow{d}\\
		\mathbb N^n \arrow{r}{e_i \mapsto x_i} & k[x_1\dots x_n,T]/(x_1\dots x_n-T).
	\end{tikzcd}
	 \end{center}
	% where the left vertical map is the diagonal map.
	 	\item Retain the notation of Example \ref{examplelogstructure}(4) and take the fiber $X^{log}_0\rightarrow 0^{log}$ at $0$ of the morphism $X^{log}\rightarrow C^{log}$. Then the situation is described by the commutative diagram   	   
	 	\begin{center}
			 \begin{tikzcd}
		\mathbb N \arrow{r}{1 \mapsto 0} \arrow{d}{1 \mapsto \sum e_i}& k\arrow{d}\\
		\mathbb N^n \arrow{r}{e_i \mapsto x_i} & k[x_1\dots x_n]/(x_1\dots x_n).
	\end{tikzcd}
	 \end{center}
	% where the left vertical map is the diagonal map.
	 \item Retain the notation of Example \ref{examplelogstructure}(5), set $J=\{1,\dots, r\}$ and endow the open strata $X^0_J:=\Spec(k[x^{\pm 1}_{r+1}\dots x^{\pm 1}_n])\subseteq X_0$, with the log structure induced by the pre-log structure $\mathbb N^r\rightarrow k[x^{\pm 1}_{r+1}\dots x^{\pm 1}_n]$ sending $e_i$ to $0$. We write $X_J^{0,log}$ for the corresponding log-variety. The natural morphism $X^{0,log}_J\rightarrow 0^{log}$ is described by the commutative diagram
   	   	   \begin{center}
			 \begin{tikzcd}
		\mathbb N \arrow{r}{1 \mapsto 0} \arrow{d}{1 \mapsto \sum e_i}& k\arrow{d}\\
		\mathbb N^{r} \arrow{r}{e_i \mapsto 0} & k[x^{\pm 1}_{r+1}\dots x^{\pm 1}_n].
	\end{tikzcd}
	 \end{center}
	
\end{enumerate}
\end{example}
\subsubsection{The case of a semistable family}\label{subsec: semistablecaselog}
More generally, let $C$ be a smooth $k$-curve and $0\in C(k)$. Let $X\rightarrow C$ be a projective family, smooth over $U_C:=C\setminus \{0\}$ and having strictly semistable reduction at $0$. 

The natural inclusion $\mathcal O_C\cap j_*\mathcal O_{U_C}^*\rightarrow \mathcal O_C$ is induces a log structure on $C$, called the log-structure $C^{\log}$ associated to the divisor $0$ in $C$. Locally around $0$, in the \'etale topology, the logarithmic structure is given by the one in Example \ref{examplelogstructure}(3); see \cite[Example 1.5]{Kato89}.

Let $U_X:=X\setminus X_0\subseteq X$. The natural inclusion $\mathcal O_X\cap j_*\mathcal O_{U_X}^*\rightarrow \mathcal O_X$ is a pre-log structure on $X$, which induces a log-strucure $X^{log}$ giving rise to a morphism of log-variety $X^{\log}\rightarrow C^{\log}$, which, locally around $0$, is isomorphic to the one given in Example \ref{examplelogstructure}(5). The restriction of such log structure to $X_0$ induces a canonical log-structure $X^{\log}_0$ on $X_0$ and a morphism $X^{\log}_0\rightarrow \Spec(k)^{\log}$ (see \cite[Example 3.7 (2)]{Kato89}), which, locally around $0$, is isomorphic to the one given in Example \ref{examplelogstructure}(6). 

Endowing the strata $X^0_J\subseteq X_0$ with the induced log-structure, we get a strict closed immersion of log-schemes $X^{0,log}_J\subseteq X_0^{log}$, which is locally isomorphic to the one in example \ref{examplelogstructure}(7). The situation is summarized in the following commutative diagram:
	\begin{center}
			 \begin{tikzcd}[row sep= huge]
			 X_{J}^{0,log}\arrow{d}\arrow{rrr}{\pi} & & &X_{J}^{0}\arrow{d} \\
			 
		X_0^{log}\arrow[bend left]{rrr}{\pi}\arrow{r}\arrow{d}{f^{log}_0}	& X^{log}\arrow{r}{\pi}\arrow{d}{f^{log}}&X^{}\arrow{d}{f}& X_0\arrow{l}\arrow{d}{f_0}\\
	
					\Spec(k)^{log}\arrow[bend right]{rrr}{\pi}\arrow{r}{0}& 	  C^{log}\arrow{r}{\pi}  & C & \Spec(k)\arrow{l}{0}
	\end{tikzcd}
		 \end{center} 
 \subsection{Analytification}\label{subsec: analytification}
Assume now that $k=\mathbb C$ and $(X,M)$ is a log $\C$-variety. We recall how to associate to $(X,M)$ a topological space, the Kato-Nakayama space of $(X,M)$, which is a topological incarnation of $(X,M)$ extending the usual analytification for $\mathbb C$-varieties. We refer the reader to \cite{KatNak99} for the original construction and to \cite[Section 5.5]{illusielog} for a nice introduction.
\subsubsection{Definitions and examples}
More precisely, there is a natural functorial way to associate to $(X,M)$ a $C^{\infty}$-manifold with corner  $(X,M)^{an}$, called the Kato-Nakayama space of $(X,M)$. Concretely, consider the pair $(\Spec(\C),\pi)$ where $\pi:\R_{\geq 0}\times S^1\rightarrow \C$ is the polar coordinates defined by $\pi(x,y)=xy$. Then, by definition 
$$(X,M)^{an}:=\Hom((\Spec(\C),\pi),(X,M))$$ 
endowed with a natural topology. We will often write $(X,M)(\mathbb C)$ for $(X,M)^{an}$. 

If $M$ is the trivial log-structure, then $(X,M)^{an}=X(\mathbb C)$ is the usual analytification. In general, for every log-variety $(X,M)$, the natural morphism $(X,M)\rightarrow X$ induces a proper morphism $(X,M)^{an}\rightarrow X^{an}$, whose fiber over a point $x\in X^{an}$ is $(S^1)^{\rank(\overline M^{gp}_x)}$. 
Here some examples which are relevant for this paper.
\begin{example}\label{exampleanilytification}
\begin{enumerate}
	\item[]
	\item As already mentioned, if $X$ has the trivial log structure then $X^{an}=X(\mathbb C)$.
	\item In the notation of Example \ref{examplelogstructure}(2), $(Spec(\mathbb C)^{\log})^{an}=S^1$, since the only morphisms $h:\mathbb N\rightarrow \R_{\geq 0} \times S^1$ making the diagram
	    \begin{center}
			 \begin{tikzcd}
		\mathbb N \arrow{r}{h} \arrow{d}{1 \mapsto 0}& \R_{\geq 0} \times S^1 \arrow{d}{(x,y) \mapsto  xy}\\
		\mathbb \C \arrow{r}{id} & \C.
	\end{tikzcd}
	 \end{center}
	 commutative are those sending $1$ to $(0,y)$ for $y\in S^1$.
	\item In the notation of Example \ref{examplelogstructure}(3), then $C^{log,an}\simeq \R_{\geq 0}\times S^1$ is the real blow-up of $\mathbb C$ at $0$, as one readily sees by looking at the set of commutative diagrams
		    \begin{center}
			 \begin{tikzcd}
		\mathbb N \arrow{r} \arrow{d}{1\mapsto T}& \R_{\geq 0} \times S^1 \arrow{d}{(x,y) \mapsto  xy}\\
		\mathbb \C[T] \arrow[two heads]{r}{t\mapsto a_i} & \C.
	\end{tikzcd}
	 \end{center}
	 Moreover, the natural map $C^{log,an}\simeq \R_{\geq 0}\times S^1\rightarrow \C\simeq C^{an}$ is the usual real-blow-up map (the polar coordinates).
\item Retain the notation of Example \ref{examplelogstructure}(6). 
Then $X^{0,\log}_{J}(\mathbb C)\simeq (S^1)^r\times (\C^*)^{n-r}$, as one readily sees by looking at the set of commutative diagrams
		    \begin{center}
			 \begin{tikzcd}
		\mathbb N^{r} \arrow{r} \arrow{d}{e_i\mapsto 0}& \R_{\geq 0} \times S^1 \arrow{d}{(x,y) \mapsto  xy}\\
		\mathbb \C[x_{r+1}^{\pm 1}\dots x_n^{\pm 1}] \arrow[two heads]{r}{t\mapsto a_i} & \C.
	\end{tikzcd}
	 \end{center}
  and the map $X^{0,\log}_{J}(\mathbb C)\rightarrow X^{an}\simeq (\C^*)^{n-r}$ is the natural projection. 
  
  The natural map $X^{0,\log}_{J}\rightarrow \Spec(\C)^{\log}$ induces a map $X_J^{0,\log}(\mathbb C)\rightarrow S^1$, which locally sends $((y_1\dots y_r),(a_{r+1}\dots a_{n}))$ to $\prod_i y_i$. Indeed, $((y_1\dots y_r),(a_{r+1}\dots a_{n}))$ corresponds to the commutative diagram  
  \begin{center}
			 \begin{tikzcd}
		\mathbb N^{r} \arrow{r}{e_i\mapsto (0,y_i)} \arrow{d}{e_i\mapsto 0}& \R_{\geq 0} \times S^1 \arrow{d}{(x,y) \mapsto  xy}\\
		\mathbb \C[x_{r+1}^{\pm 1}\dots x_n^{\pm n}] \arrow[two heads]{r}{t\mapsto a_i} & \C.
	\end{tikzcd}
	 \end{center}
	 and its image in $S^1$ it is obtained by precomposing it with the diagram,
	   \begin{center}
			 \begin{tikzcd}
			 \mathbb N \arrow{r}{1 \mapsto \sum e_i}\arrow{d}{1\mapsto 0}& \mathbb N^{r}\arrow{d}{e_i\mapsto 0} \\
			 \C\arrow{r} & \C[x_{r+1}^{\pm 1}\dots x_n^{\pm n}]\\
	\end{tikzcd}
	 \end{center}
  in which the bottom horizontal map is the canonical inclusion. Finally, one obtains the diagram
    		    \begin{center}
			 \begin{tikzcd}[column sep= huge]
		\mathbb N \arrow{r}{e_i\mapsto (0,\prod_i y_i)} \arrow{d}& \R_{\geq 0} \times S^1 \arrow{d}{(x,y) \mapsto  xy}\\
		\mathbb \C \arrow[equal]{r} & \C,
	\end{tikzcd}
	 \end{center}
  which corresponds to the point $\prod y_i\in S^1$. In particular, the fiber $X^{0,\log}_{J,1}(\mathbb C)$ at $1$ of this morphism is the subset of $X^{0,\log}_{J}(\mathbb C)$, in which the product of the first $r$ coordinates is $1$.

\end{enumerate}
\end{example}
\subsubsection{The case of a semistable degeneration}\label{subsec: semistablecasean}
Assume from now that $X\rightarrow C$ is a real strictly semistable degeneration as in Section \ref{subsec: semistablecaselog} and endow it with the logarithmic structure as in Section \ref{subsec: semistablecaselog}. Choose an orientation $(U,\psi)$  around $0$, see Section (\ref{sec : orientation}). Then, by Example \ref{exampleanilytification}(2), the orientation induces an isomorphism
$\pi^{-1}(U)\simeq [0,1)\times S^1$, where $\pi:C^{log,an}\rightarrow C^{an}$ is the natural map. Hence $C^{log,an}$ can be identified with the real oriented blow-up of $C^{an}$ in $0$ (see e.g. \cite[(1.2.3)]{nakayamaogus}).

Then, we can analytify the commutative diagram
	\begin{center}
			 \begin{tikzcd}
			 X_{J}^{0,log}\arrow{d}\arrow{r}{\pi}& X^{0}_J\arrow{d}\\
			 X_0^{log}\arrow{d}\arrow{r}{\pi}& X_0\arrow{d}\\
			 0^{log}\arrow{r}& 0
	\end{tikzcd}
		 \end{center} 
to get morphisms $X_{J}^{0,log}(\mathbb C)\rightarrow X_0^{log}(\mathbb C)\rightarrow 0^{log,an}\simeq S^1$ with fibers at $1$ given by $X_{J,1}^{log}(\mathbb C)\rightarrow X_{0,1}^{log}(\mathbb C)$ and a commutative cartesian diagram
	\begin{center}
			 \begin{tikzcd}
		X_{J,1}^{0,log}(\mathbb C)\arrow{r}\arrow{d}& X_{J}^{0,log} (\C)\arrow{d}\arrow{r}{\pi}& X^{0}_J(\mathbb C)\arrow{d}\\
		X_{0,1}^{log}(\mathbb C)\arrow{d}\arrow{r}&	 X_0^{log}(\mathbb C)\arrow{d}\arrow{r}{\pi}& X_0(\mathbb C)\\
		\{1\}\arrow{r}&	 0^{\log}(\mathbb C)=S^1.
	\end{tikzcd}
		 \end{center} 
The main result that we are interested in and which motivates the use of logarithmic geometry, is the following.
\begin{theorem}\cite[Theorem 0.3]{nakayamaogus}\label{ogusnakayama}
The topological space $X^{log}_{0,1}(\mathbb C)$	 is homeomorphic to every smooth fiber of the morphism $X\rightarrow C$.
\end{theorem}
Hence, in order to study the Betti numbers of the complex general fiber of $X\rightarrow C$, we can study the topology of $X^{log}_{0,1}(\mathbb C)$ via the morphism $X^{log}_{0,1}(\mathbb C)\rightarrow X_0(\mathbb C)$.
\subsection{Stratifications}\label{subsec: stratification}
\numberwithin{equation}{subsection} 
We now go further and we study $X_{0,1}^{log}(\mathbb C)$, by stratifying it in a way which is compatible with the logarithmic structure on $X_{0}(\mathbb C)$. Retain the notation of Section \ref{subsec: semistablecasean}. 

Let $\mathfrak Z:=\{X^0_\Delta\}_{\Delta \in P}$ be an algebraic refinement of $\mathfrak I$. Let $I_{\Delta}$ be the subset $I_\Delta\subseteq I$ such that $X^0_\Delta\subseteq X^0_{I_\Delta}$ and set $\vert \Delta\vert:=\vert I_\Delta\vert$.
 We endow $X^0_\Delta$ with the structure of a log-variety $X^{0,log}_\Delta$ by restricting the log structure of $X_0$ to $X^0_\Delta$.
 The most important properties of $X^{0,log}_\Delta(\mathbb C)$ and $X^{0,log}_{\Delta,1}(\mathbb C)$, which follows directly from the local description given in Example \ref{exampleanilytification} are summarized in the following remark.
 \begin{remark}
 \begin{enumerate}\label{keyproperty}
 \item[]
 	\item Example \ref{exampleanilytification}(7) shows that $X^{0,log}_\Delta(\mathbb C)\rightarrow X^0_\Delta(\mathbb C)$ is a locally trivial fibration with fibers $(S^1)^{\vert \Delta\vert}$. In particular $X^{0,log}_\Delta(\mathbb C)$ is a smooth $C^{\infty}$-manifold of dimension $2\Dim(X_0)+\vert\Delta \vert$.
 	\item By Example \ref{exampleanilytification}(4), the morphism  $X^{0}_{\Delta,1}(\mathbb C)\rightarrow X^{0}_{\Delta}(\mathbb C)$ is a proper locally trivial fibration with fiber $(S^1)^{\vert \Delta\vert-1}$. In particular, $X^{0,log}_\Delta(\mathbb C)$ is a smooth $C^{\infty}$-manifold of dimension $2\Dim(X_0)+\vert\Delta \vert-1$.
 \end{enumerate}
 \end{remark}
Since $\{X^{0}_{\Delta}(\mathbb C)\}$ is a stratification of $X_0$ 
and $\overline{X_{\Delta}^{0,\log}(\mathbb C)}=X_{\Delta}^{\log}(\mathbb C)\subseteq X^{log}_{0}(\mathbb C)$\footnote{This can be checked locally, hence we can assume that $X_0^{\log}=Spec(\mathbb C[x_1,\dots , x_n]/(x_1\dots x_n))$ endowed with the logarithmic structure of Example \ref{examplelogstructure}(5). Since in this case  $X_0^{\log}(\mathbb C)$ is closed inside $\mathbb A^{n,log}(\mathbb C)\simeq (\R_{\geq 0}\times S^1)^n$, we can assume that $X_0^{log}$ is $\mathbb A^{n}=\prod \mathbb A^1$ with the product log-structure of example \ref{examplelogstructure}(2). In this case the result follows from a computation using that the map $(\R_{\geq 0}\times S^1)^n\simeq \mathbb A^{n,\log}(\mathbb C)\rightarrow \mathbb A^{n}(\mathbb C)=\mathbb C^n$ is the polar coordinates componentwise.} , we see that $\{X^{0,log}_{\Delta}(\mathbb C)\}$ is a stratification of $X^{log}_{0}(\mathbb C)$.
Since $X^{log}_{0,1}(\mathbb C)\subseteq X^{log}_{0}(\mathbb C)$ is a closed subset, also $\{X^{0,log}_{\Delta,1}(\mathbb C)\}$ is a stratification of $X^{log}_{0,1}(\mathbb C)$.
\subsection{Construction of the complex and statement of the main theorem}\label{subsec: construction}
Retain the notation of Section \ref{subsec: stratification} and let $A$ be a coefficient ring. 
There is a natural spectral sequence (see e.g. \cite{Pete16} and \cite[(3) pag. 2527]{Pete17})
\begin{equation}\label{eq : spectral complex}
	^{\C}E_1^{a,b}:=\bigoplus_{\Dim(X^0_\Delta)=a} H_c^{a+b}(X^{0,log}_{\Delta,1}(\mathbb C),A)\Rightarrow H_c^{a+b}(X_{0,1}^{log}(\mathbb C),A)\big (\simeq H_c^{a+b}(X_t(\mathbb C),A)\big ),
\end{equation}
where $t\in C(\mathbb C)$ is any member different from $0$ and the last isomorphism follows from Theorem \ref{ogusnakayama}.

Set 

$$C_{q,\mathfrak Z,A}^{\bullet}:=\text{} ^{\C}E_1^{\bullet,q}:$$
$$ \underset{\Dim(X^0_\Delta)=p-1}\bigoplus H_c^{p+q-1}(X^{0,log}_{\Delta,1}(\mathbb C),A)\rightarrow \underset{\Dim(X^0_\Delta)=p}\bigoplus H_c^{p+q}(X^{0,log}_{\Delta,1}(\mathbb C),A)\rightarrow
			 \underset{\Dim(X^0_\Delta)=p+1}\bigoplus H_c^{p+q+1}(X^{0,log}_{\Delta,1}(\mathbb C),A).$$
			 
Now Theorem \ref{main} can be stated more precisely as follows.
\begin{theorem}\label{mainpreciso}
\begin{enumerate}
\item[]	
	\item With the notation of Theorem \ref{main}, assume that 
 	\begin{enumerate}
\item[(a)] $H^i(X^0_\Delta(\mathbb R),\mathbb Z/2\Z)=0$ for all $i\geq 1$ and $X^0_\Delta\in \mathfrak Z$;
 \item[(b)]  $X^0_\Delta$ is a maximal variety, for all $X^0_\Delta\in \mathfrak Z$.
 	\end{enumerate}
 	
 	Then, for every $t\in C(\mathbb R)$ near $0$ one has:
 	$$b_p(X_t(\mathbb R)),\Z/2\Z)\leq \sum_{0\leq q\leq n}\Dim(H^p(C_{q,\mathfrak Z,\mathbb Z/2\Z}^{\bullet}).$$
\item If in addition to (a) and (b), assume that
 	\begin{enumerate}
 		\item[(c)] the mixed Hodge structure on  $H^i(X^0_\Delta,\Q)$ is pure of type $(i,i)$ and $H^i(X^0_\Delta,\Z)$ is torsion free,
 	 	 	\end{enumerate}  
 	 	 	
 	Then, for every $t\in C(\mathbb R)$ near $0$ one has:
 	
	$$(i)\quad \Dim(H^q(C_{p,\mathfrak Z,\Z}^{\bullet}\otimes \Q))=h^{p,q}(X_t)\quad \text{and} \quad (ii) \quad C_{p,\mathfrak Z,\Z}^{\bullet}\otimes \Z/2\Z \simeq C_{p,\mathfrak Z,\mathbb Z/2\Z}^{\bullet}.$$
\end{enumerate}
 \end{theorem}
 
 	\numberwithin{equation}{subsubsection} 

  \section{Real log varieties and equivariant cohomology} \label{sec: real_log}
  In this section we recall some facts from the theory of real logarithmic varieties. The theory appears explicitly in the literature in \cite[Sections 5-7]{Argu21} and in \cite{ArgBou21}, and it is implicit in \cite{Rau22}. We borrow part of the presentation from \cite[Sections 5-7]{Argu21}, where the author works in a slightly more general setting and we complement Arg\"uz results adding what is needed for our purposes. In Section  \ref{subsec : deflogreal}, we recall the basic definitions,  examples and we study the case of totally real semistable degenerations. In Section \ref{subsec : stratlogreal} we study stratifications of these degenerations. In Section \ref{subsec : equivariant} we recall a few fact around equivariant cohomology.
    \subsection{Real logarithmic structure}\label{subsec : deflogreal}
    \subsubsection{Definitions and examples}
Assume now that $k=\mathbb R$. To give a log variety over $\R$ is equivalent to give a log variety $(X,M)$ over $\C$ endowed with a pair $(\sigma, \sigma_M)$ where $\sigma: X(\mathbb C)\rightarrow X(\mathbb C)$ is an anti-holomorphic involution and $\sigma_M:M\rightarrow M$ is an involution of monoids making the diagram  \begin{center}
			 \begin{tikzcd}
		M\arrow{d}\arrow{r}{\sigma_M} & M\arrow{d}\\
		\mathcal O_X \arrow{r}{\sigma} & \mathcal O_X . 
	\end{tikzcd}
		 \end{center}
		  commutative, where, by abuse of notation, we denote by $\sigma$ also the involution on $\mathcal O_X$ induced by that on $X(\mathbb C)$.

		 We endow $\R_{\geq 0}\times S^1$ with the involution $\tau$ given by $\tau(x,r)=(x,\overline r)$, so that the polar coordinates morphism $\R_{\geq 0}\times S^1\rightarrow \mathbb C$ is equivariant once $\mathbb C$ is endowed with the standard conjugation. Then $\tau$ and $\sigma$ induce a natural involution on $(X,M)(\mathbb C):=\Hom((\Spec(\C),\pi),(X,M))$ and  the real locus $(X,M)(\mathbb R)\subseteq (X,M)(\mathbb C)$ of $(X,M)$ is, by definition, the fixed locus of this involution endowed with its subspace topology. 
		 
		 If $(X,M)\rightarrow (Y,N)$ is a morphism of real log varieties, then the induced morphism $(X,M)(\mathbb C)\rightarrow (Y,N)(\mathbb C)$ is equivariant and hence it induces a morphism $(X,M)(\mathbb R)\rightarrow (Y,N)(\mathbb R)$. In particular, there is a natural morphism $(X,M)(\mathbb R)\rightarrow X(\mathbb R)$. We now give the examples that are relevant for our paper.  
		 \begin{example}\label{examplerealanalityfication}
		 \begin{enumerate}
		 	\item[]
		 	\item If $X$ is endowed with the trivial log-structure, then $(X,M)(\mathbb R)\simeq X(\mathbb R)$ is just the set of real points of $X$.
		 	\item Let $X=\Spec(\mathbb R)$ endowed with the canonical log-structure. Then the induced action on $(X(\mathbb C),\mathbb N)$ is trivial on both $X(\mathbb C)=\{\text{pt}\}$ and $\mathbb N$. This induces on the analytification $X^{\log}(\mathbb C)\simeq S^1$ the standard conjugation, so that $X^{\log}(\mathbb R)=\{+1,-1\}$.
		 	\item Assume that $C=\Spec(\mathbb R[T])$, endowed with the log structure from Example \ref{examplelogstructure}(3). Then, the involution on $(C(\mathbb C),\mathbb N)$ is the standard complex conjugation on $C(\mathbb C)\simeq \mathbb C$ and the identity on $\mathbb N$. This induces on the analytification $C^{\log}(\mathbb C)\simeq \R_{\geq 0}\times S^1$ (see Example \ref{exampleanilytification}(3)) the identity of the first factor and the standard complex conjugation on the second one, so that $C^{log}(\mathbb R)\simeq \R_{\geq 0}\times \{\pm 1\}$. The natural map $C^{log}(\mathbb R)\simeq \R_{\geq 0}\times \{\pm 1\}\rightarrow C(\mathbb R)\simeq \R$ sends $(x,1)$ to $x$ and $(x,-1)$ to $-x$. 
		 	\item Assume that $X^0_J=\Spec(\mathbb R[x^{\pm 1}_{r+1} \dots x^{\pm 1}_{n}])$, endowed with the log structure $\mathbb N^r\rightarrow \mathbb R[x^{\pm 1}_{r+1} \dots x^{\pm 1}_{n}]$ from Example \ref{examplelogstructure}(6) sending all non zero elements to $0$. Then the involution on $(X^0_J(\mathbb C),\mathbb N^r)$ is the standard complex conjugation on $X^0_J(\mathbb C) \simeq (\C^{*})^{n-r}$ and the trivial action on $\mathbb N^r$. Hence the induced action on $X_J^{0,\log}(\mathbb C)\simeq (S^1)^{r}\times (\C^{*})^{n-r}$ (see Example \ref{exampleanilytification}(3)) is the standard complex conjugation acting component-wise and without any permutation of the factors, so that $X_J^{0,\log}(\mathbb R)\simeq \{\pm 1\}^{r}\times (\R^{*})^{n-r}$ and the map $X_J^{0,\log}(\mathbb R) \simeq\{\pm 1\}^{r}\times (\R^{*})^{n-r}\rightarrow X(\R)\simeq (\R^{*})^{n-r}$ is the natural projection.  The natural map $X^{0,\log}_J(\R)\simeq \{\pm 1\}^{r}\times (\R^{*})^{n-r} \rightarrow \Spec(\mathbb R)^{\log}\simeq \{\pm 1\}$ sends $((y_1\dots y_r),(a_{r+1}\dots a_{r}))$ to $\prod_i y_i$, as follows from the computation done in Example \ref{exampleanilytification}(4). In particular, the fiber $X_{J,1}^{0,\log}(\mathbb R)$ at $1$ of this morphism identifies with the subset of  $X^{0,\log}_J(\mathbb R)\simeq \{\pm 1\}^{r}\times (\R^{*})^{n-r}$ in which the product of the first $r$ components is $1$. Hence  $X^{0,\log}_{J,1}(\mathbb R)\rightarrow X^{0}_J (\mathbb R)$ is a topological cover of degree $2^{\vert J\vert -1}$; see also \cite[Proposition 8.5]{Argu21} and \cite[Theorem 1.3]{Rau22}. 
		 \end{enumerate}
		 \end{example}
		 In all the examples of Example \ref{examplerealanalityfication}, the action of the involution on the sheaf of monoids is trivial. Even if this is the only case we will need (since the degenerations we are interested in are totally real), for completeness, we give now an example in which the action is non-trivial, which comes from a real semistable degeneration which in not totally real. The reader only interested in Theorem \ref{main} can skip the end of this section.
		 \begin{example}
		 Let $X^0_J$ be $\Spec(\C[x_{r+1}^{\pm 1}\dots x_n^{\pm 1}])$ and endow it with the logarithmic structure associated to the morphism $\mathbb N^{r}\rightarrow \Spec(\C[x_{r+1}^{\pm 1}\dots x_n^{\pm 1}])$ sending all the non zero element to $0$. We consider the standard conjugation action on $X^0_J(\mathbb C)$ but now we consider $\mathbb N^{r}$ endowed with the involution switching $e_1$ and $e_2$ and acting trivially on the other $e_i$. This corresponds to the restriction to some geometrically irreducible strata of  the log structure coming from a degeneration whose special fiber has the form $\Spec(\mathbb R[x_1\dots x_n]/((x^{2}_1+x_2^{2})x_3\dots x_n)$, which is not totally real.

		 Then the induced involution on $X^0_J(\mathbb C)\simeq (S^1)^{r}\times (\C^{*})^{n-r}$ sends $((y_1,y_2,y_3\dots y_r),(a_{r+1}\dots a_n)) $ to
		 $$((\overline{y_2},\overline{y_1},\overline{y_3}\dots \overline{y_r}),(\overline{a_{r+1}}\dots \overline{a_n})),$$
		 so that $X^{0,\log}_J(\mathbb R)\subseteq X^{0,\log}_J(\mathbb C)$ identifies with the subset made by the elements of the form 
		 $$((y_1,\overline{y_1},\pm 1,\dots, \pm 1),(a_{r+1}\dots a_n))$$ 
		 with $y_1\in S^1$ and $a_{i}\in \R^*$. Hence $X^{0,\log}_J(\mathbb R)$ identifies with $S^1\times \{\pm 1\}^{r-2}\times (\R^{*})^{n-r}$.

		 The natural map $S^1\times \{\pm 1\}^{r-2}\times (\R^{*})^{n-r}\simeq X^{0,log}_J(\mathbb R)\rightarrow X^0_J(\mathbb R)\simeq (\R^{*})^{n-r}$ is the natural projection, which is no longer a topological cover of degree $2^{r}$ but a product of a $2^{r-2}$ topological cover and an $S^1$-bundle  over $X^0_J(\mathbb R)$. The fiber at $1$ of the other natural map $S^1\times \{\pm 1\}^{r-2}\times (\R^{*})^{n-r} \simeq X^{0,log}_J(\mathbb R) \rightarrow \Spec(\mathbb R)^{\log}\simeq \{\pm 1\}$ is given (since $y_1\overline{y_1}=1$) by the subset of $X^{0,log}_J(\mathbb R)$ made by elements $((y_1,\overline{y_1},y_3\dots y_r),(a_{r+1}\dots a_n))$ 
		 for which $\prod_{i\geq 3}(y_i)=1$, hence it is a product of an $S^1$-bundle and a topological cover of degree $2^{r-3}$ over $X^0_J(\mathbb R)$. 
		 \end{example}
	    \subsubsection{The case of a strictly semistable degeneration}\label{subsec : semistablecase real}
Assume from now that $X\rightarrow C$ is a totally real strictly semistable degeneration and endow it with the logarithmic structures as in Section \ref{subsec: semistablecaselog}. Then we have a morphisms of real log-varieties $X^{log}\rightarrow C^{log}$ and  $X_0^{log}\rightarrow 0^{log}$. Locally around $0$, the log structure are isomorphic to one given in Example \ref{examplerealanalityfication}. Fix an orientation $(U,\phi)$ around $0$; see Section \ref{sec : orientation}.

 	 Recall the commutative diagram from Section \ref{subsec: semistablecasean} 
 		\begin{center}
			 \begin{tikzcd}
		X_{J,1}^{0,log}(\mathbb C)\arrow{r}\arrow{d}&	 X_{J}^{0,log}(\mathbb C)\arrow{d}\arrow{r}{\pi}& X^{0}_J(\mathbb C)\arrow{d}\\
		X_{0,1}^{log}(\mathbb C)\arrow{d}\arrow{r}&	 X_0^{log}(\mathbb C)\arrow{d}\arrow{r}{\pi}& X_0(\mathbb C)\\
		\{1\}\arrow{r}&	 0^{\log}(\mathbb C)=S^1.
	\end{tikzcd}
		 \end{center} 
The morphisms on the right part of the diagram are compatible with the involutions acting on the various topological spaces appearing. Hence the involutions acting on $X^{0,log}_J(\mathbb C)$ and $X^{log}_0(\mathbb C)$ restrict to involutions on $X^{0,log}_{J,1}(\mathbb C)$ and $X_{0,1}^{log}(\mathbb C)$ respectively, making the whole diagram equivariant.

Taking fixed points we get a commutative diagram
 		\begin{center}
			 \begin{tikzcd}
		X_{J,1}^{0,log}(\mathbb R)\arrow{r}\arrow{d}&	 X_{J}^{0,log}(\mathbb R)\arrow{d}\arrow{r}{\pi}& X^{0}_J(\mathbb R)\arrow{d}\\
		X_{0,1}^{log}(\mathbb R)\arrow{d}\arrow{r}&	 X_0^{log}(\mathbb R)\arrow{d}\arrow{r}{\pi}& X_0(\mathbb R)\\
		\{1\}\arrow{r}&	 0^{\log}(\mathbb R)=\{\pm 1\},
	\end{tikzcd}
		 \end{center} 
 	 	where $X_{J,1}^{0,log}(\mathbb R)$ and $X_{0,1}^{log}(\mathbb R)$ are both the fiber at $1$ of the morphisms 
 	 	$X_{J}^{0,log}(\mathbb R)\rightarrow X^{log}_0 (\mathbb R)\rightarrow \{\pm 1\}$ 
 	 	and the fixed points of the involutions acting respectively on $X^{0,log}_{J,1}(\mathbb C)$ and $X_1^{0,log}(\mathbb C)$.
 The main result that we are interested in and which motivates the use of real logarithmic geometry, is the real analogue of \cite{nakayamaogus}, which follows from the results in \cite{Argu21} or in \cite{Rau22}. 
\begin{theorem}[\cite{Argu21}, \cite{Rau22}]
\label{katonakayamareal}
The topological space $X^{log}_{0,1}(\mathbb R)$ is homeomorphic to every positive real smooth fiber near $0$ of the morphism $X\rightarrow C$.
\end{theorem}
\proof
This is essentially the content of \cite[Proposition 7.4]{Argu21} or the one of \cite[Theorem 1.1]{Rau22}, which say
\footnote{We explain how to deduce this from the results there.
\begin{itemize}
	\item Even if \cite[Proposition 7.4]{Argu21} is stated in a very particular case, the proof extends to nice (i.e. exact) proper and logarithmic smooth morphism, since it (explicitly) mimics the proof of Theorem \ref{katonakayamareal} in \cite{nakayamaogus};
	\item \cite[Theorem 1.1]{Rau22} is stated for morphisms around a disk. Since an orientation around $0$ induces an equivariant isomorphism between a local neighbourhood of $0$ with a disk and of a local neighbourhood of $0$ in $C^{log}(\mathbb R)$ with the real oriented blow-up, Rau's results apply to this situation.
\end{itemize}}
that $X^{log}(\mathbb R)\rightarrow C^{log}(\mathbb R)$ is a locally trivial fibration. 
Therefore all smooth real fibers in the same connected component of $C^{log}(\mathbb R)$ are homeomorphic.
Since the orientation induces an equivariant isomorphism 
$$\pi^{-1}(U)\simeq [0,1)\times S^1 \quad\text{so that}\quad \pi_{\mathbb R}^{-1}(U(\mathbb R))\simeq [0,1)\times \{\pm 1\},$$
and the fibers of $X^{log}(\mathbb R)\rightarrow C^{log}(\mathbb R)$ are homeomorphic to the ones of $X(\mathbb R)\rightarrow C(\mathbb R)$ outside $0$, the positive fibers of $X(\mathbb R)\rightarrow C(\mathbb R)$ near $0$ are all homeomorphic to $X^{log}_{0,1}(\mathbb R)$.
\endproof
Hence, in order to study the Betti numbers of a positive smooth real general fiber of $X(\mathbb R)\rightarrow C(\mathbb R)$, we can study the topology of $X^{log}_{0,1}(\mathbb R)$ via the morphism $X^{log}_{0,1}(\mathbb R)\rightarrow X_0(\mathbb R)$.	 	
 	 	
	\subsection{Real stratifications}\label{subsec : stratlogreal}
	\numberwithin{equation}{subsection} 

    Let $\mathfrak Z:=\{X^0_\Delta\}_{\Delta \in P}$ be an algebraic refinement of $\mathfrak I$ made by smooth real algebraic varieties and let us retain the notation of Section \ref{subsec: stratification}. 
 We endow $X^0_\Delta$ with the structure of a real log-variety $X^{0,log}_\Delta$ by restricting the real log structure of $X_0$ to $X^0_\Delta$.
 The most important properties of $X^{0,log}_\Delta(\mathbb R)$ and $X^{0,log}_{\Delta,1}(\mathbb R)$, which follows directly from the local description given in Example \ref{exampleanilytification}, are summarized in the following remark.
 \begin{remark}
 \begin{enumerate}\label{keypropertyreal}
 \item[]
 	\item Example \ref{exampleanilytification}(4) shows that $X^{0,log}_\Delta(\mathbb R)\rightarrow X^0_\Delta(\mathbb R)$ is a topological cover of degree $2^{\vert \Delta\vert}$. In particular $X^{0,log}_\Delta(\mathbb R)$ is a smooth $C^{\infty}$-manifold of dimension $\Dim(X_0)$.
 	\item Example \ref{exampleanilytification}(4) shows that the morphism  $X^{0}_{\Delta,1}(\mathbb R)\rightarrow X^{0}_{\Delta}(\mathbb R)$  is a topological cover of degree $2^{\vert \Delta\vert-1}$. In particular $X^{0,log}_{\Delta,1}(\mathbb R)$ is a smooth $C^{\infty}$-manifold of dimension $\Dim(X_0)$.
 \end{enumerate}
 \end{remark}
By Section \ref{subsec: stratification}, one has that $\{X^{0,\log}_{\Delta,1}(\mathbb C)\}$ is a stratification of $X^{log}_{0,1}(\mathbb C)$. Since $X^{log}_{0,1}(\mathbb R)\subseteq X^{log}_{0,1}(\mathbb C)$ is closed, also $\{X^{0,\log}_{\Delta,1}(\mathbb R)\}$ is a stratification of $X^{log}_{0,1}(\mathbb R)$.

We can now compute the Betti numbers of $X^{0,log}_{\Delta,1}(\mathbb R)$ in terms of the Betti numbers of $X^0_{\Delta}(\mathbb R)$. This will be important to transfer the assumptions $(a)$ and $(b)$ in Theorem \ref{main} from $X^0_{\Delta}(\mathbb C)$ to $X^{0,log}_{\Delta,1}(\mathbb C)$.
\begin{lemma}\label{totalbettinumberreale}
	Assume that $H^1(X^0_{\Delta}(\mathbb R),\mathbb Z/2\mathbb Z)=0$. Then for every $i\geq 0$ one has
		$$b_i(X^{0,log}_{\Delta,1}(\mathbb R),A)=2^{\vert\Delta\vert-1}b_i(X^0_{\Delta}(\mathbb R),A).$$
\end{lemma}
\proof
By Remark \ref{keypropertyreal}, $X^0_{\Delta,1}(\mathbb R)\rightarrow X^0_{\Delta}(\mathbb R)$ is a topological cover of degree $2^{\vert\Delta\vert-1}$. Since $H^1(X^0_{\Delta}(\mathbb R),\mathbb Z/2\mathbb Z)=0$, this cover is the trivial one. Hence $X^0_{\Delta,1}(\mathbb R)$ is the disjoint union of $2^{\vert\Delta\vert-1}$ copies of $X^0_{\Delta}$, and the lemma follows.
\endproof
\subsection{Equivariant cohomology}\label{subsec : equivariant}
\numberwithin{equation}{subsubsection} 

In this Section we recall a few generalities on equivariant cohomology as presented in \cite[Chapter IV]{borel}. This is an important tool in the proof of Theorem \ref{mainpreciso} to compare complex and real Betti numbers. Since \cite[Chapter IV]{borel} deals with the case of cohomology and we need the compact cohomology case, we give some details.

In this section we consider $G:=\mathbb Z/2\Z$ as an abelian group and we let  $X$ be a locally compact Hausdorff  topological space on which $G$ acts.  We write $X^{G}$ for the set of fixed points for the action of $G$ and $X/G$ for the quotient of $X$. We assume that $H_c^i(X,\mathbb Z/2\Z)$ is finite dimensional for all $i\geq 0$ and that there exists some $n_X\in \mathbb N$ such that $H_c^i(X,\mathbb Z/2\Z)=0$ for all $i>n_X$.  
  \subsubsection{Group cohomology}\label{subsec : groupcoho}
    Let $M$ be a $\Z/2\Z$-vector space on which $G$ acts. Write $H^i(G,M)$ for its $i^{th}$-cohomology group. If $M$ is finite, then $H^i(G,M)$ is finite, there is an inequality $\Dim (H^i(G,M)) \leq \Dim ( M)$  and if the action of $G$ on $M$ is trivial then $H^i(G,M)=M$ (see e.g. \cite[IV, 2.1, Pag. 50]{borel}).
 \subsubsection{Classifying space}\label{subsec : class}
 We fix an auxiliary integer $N\gg 0$, e.g. $N\geq n_x+1$. Let $BG:=\mathbb R\mathbb P^{N}$ be the real projective space of dimension $N$, which is the classifying space for $G$ and for dimension $N$ (\cite[Section 18]{borel2lavendetta}), in the sense that $\pi_1(BG)\simeq G$ and for every locally constant sheaf $\mathcal F$ on $BG$ and every $x\in BG$, one has $H^i(BG,\mathcal F)=H^i(G,\mathcal F_x)$ for all $i\leq N$. Let $EG:=S^{N}$ be the unit sphere of dimension $N$, which is the universal cover $EG\rightarrow BG$ of $BG$, endowed with the natural $G$-action. 
 \subsubsection{G-equivariant cohomology with compact support} 
   Define $X_G:=(X\times EG)/G$ where $G$ acts on $X\times EG$ diagonally. Then, for every coefficient ring $A$, the equivariant cohomology of $X$ with compact support and coefficient in $A$ is defined by $H^i_{G,c}(X,A):=H_c^i(X_G,A)$.
  It is important to keep in mind the following morphisms: 
  \begin{center}
\begin{tikzcd}
X/G&X_G\arrow{r}{h}\arrow[swap]{l}{g} & BG,
	\end{tikzcd}
	 \end{center}
where  $g:X_G\rightarrow X/G$ (resp. $h:X_G\rightarrow BG$) is induced by the first (resp. the second) projection $X\times EG\rightarrow X\rightarrow X/G$ (resp. $X\times EG \rightarrow EG\rightarrow BG$).       
\subsubsection{Smith-Thom spectral sequence}
      The morphism $h:X_G\rightarrow BG$ has fibers $X$, hence the Leray spectral sequence with compact support for $h$  reads
  $$_XE_2^{a,b}:=H_c^a(BG,H_c^b(X,A))\Rightarrow H^{a+b}_{G,c}(X,A),$$
  where $H_c^b(X,A)$ is the constant local system associated to the group $H_c^b(X,A)$. Since $BG$ is compact and the classifying space for $G$ and for dimension $N$, one has $H_c^a(BG,H_c^b(X,A))=H^a(BG,H_c^b(X,A))=H^a(G,H_c^b(X,A))$; see Section \ref{subsec : class}. Therefore we get a spectral sequence
  \begin{equation}\label{smiththomspectralsequence}
  	^{ST}_{X}E_2^{a,b}:=H^a(G,H_c^b(X,A))\Rightarrow H^{a+b}_{G,c}(X,A),
  \end{equation}
  which we call the the Smith-Thom spectral sequence (with compact support). 
  \subsubsection{Maximality}\label{subsubsec maximalityequivariant}
   Recall that, by \cite[Application 3.7(b), Pag 54]{borel}, for $n_X<m<N$, the inclusion $(X^G)_G=X^G\times BG\subseteq X_G$ induces a canonical isomorphism $H^{m}_{G,c}(X,\Z/2\Z)\simeq H^m_{c}(X^G\times BG,\Z/2\Z)\simeq H_c^*(X^G,\Z/2\Z)$. Since $\Dim ( H^a(G,H_c^b(X,\Z/2\Z)) \leq \Dim ( H_c^b(X,\Z/2\Z))$, the spectral sequence (\ref{smiththomspectralsequence}) gives an inequality 
   \begin{equation}\label{smith-thom inequality}
	b_{c,*}(X^G)\leq b_{c,*}(X),
\end{equation}
where $b_{c,*}$ denotes the total Betti number with compact support. We can now recall the definition of maximality in this context.
  \begin{definition}
  	We say that $X$ is maximal (in the sense of Smith-Thom) if $b_{c,*}(X^G)=b_{c,*}(X)$. 
  \end{definition}
The Smith-Thom spectral sequence (\ref{smiththomspectralsequence}), together with the inequality $\vert H^a(G,H_c^b(X,\Z/2\Z))\vert\leq \vert H_c^b(X,\Z/2\Z)\vert$, implies that $X$ is maximal if and only if $^{ST}_{X}E_2^{a,b}$ degenerates at $E_2$ and the action of $G$ on $H_c^b(X,  \mathbb Z/ 2\mathbb Z)$ is trivial for all $b$; see e.g. \cite[4.1, Pag. 55]{borel}. 
    \subsubsection{Filtrations and maximality}\label{filtrazionest}
    Assume now that $X$ is maximal, so that the spectral sequence  degenerates at $E_2$ and the action of $G$ on $H_c^b(X, \mathbb Z/ 2\mathbb Z)$ is trivial for all $b$.
     
        Then, for every $n_X<m<N$ the spectral sequence (\ref{smiththomspectralsequence}) induces a canonical decreasing filtration 
        $$0=\text{}_mF_{m+1}\subseteq \text{}_mF_{m}\subseteq \dots \subseteq \text{}_mF_{i+1}\subseteq \text{}_mF_{i}\subseteq \dots  \subseteq \text{}_mF_{0}=H_c^*(X^G, \mathbb Z/ 2\mathbb Z)(\simeq H_{G,c}^m(X,\mathbb Z/ 2\mathbb Z))$$ such that 
    $$_mF_i/_mF_{i+1}\simeq \text{}^{ST}_{X}E_2^{i,m-i}=H^{i}(G,H_c^{m-i}(X,  \mathbb Z/ 2\mathbb Z))=H_c^{m-i}(X,  \mathbb Z/ 2\mathbb Z),$$
    where the last equality follows from the fact that the action of $G$ on $H_c^{m-i}(X, \mathbb Z/ 2\mathbb Z)$ is trivial; see Section \ref{subsec : groupcoho}.
    
  \section{Proof of Theorem \texorpdfstring{\ref{mainpreciso}}-(1)}\label{sec: equivariant_coho}
 Retain the notation and assumptions of Theorem \ref{mainpreciso}(1). Up to a choice of an orientation $(U, \psi)$ around $0 \in C(\mathbb R)$, we can assume that $X_t$, for some $t$ near $0$, is positive; see Section \ref{sec : orientation}. In order to simplify the notation, in this section, for every topological space $Z$, we write 
 $$H^i(Z):=H^i(Z,\mathbb Z/2\Z);\quad H^i_c(Z):=H_c^i(Z,\mathbb Z/2\Z);\quad b_{i,c}:=\Dim(H^i_c(Z));\quad  b_{c,*}:=\sum_i b_{i,c}.$$

  As explained in Section \ref{sec: introstrategy} the strategy is the following:
  \begin{enumerate}
  \item By using the work done in the previous sections, we construct a topological space $X_{0,1}^{log}(\mathbb C)$ endowed with an involution and with a morphism $X_{0,1}^{log}(\mathbb C)\rightarrow X_{0}(\mathbb C)$ compatible with the involutions, such that $X_{0,1}^{log}(\mathbb C)$ is homeomorphic to the general smooth fiber of the morphism $X(\mathbb C)\rightarrow C(\mathbb C)$ and the set of fixed points $X_{0,1}^{log}(\mathbb R)$ is homeomorphic to a positive real fiber near $0$;
  \item By using the stratification $\{X^0_{\Delta}\}$ and the geometry of the morphism $X_{0,1}^{log}(\mathbb C)\rightarrow X_{0}(\mathbb C)$, we stratify $X_{0,1}^{log}(\mathbb C)$ and $X_{0,1}^{log}(\mathbb R)$ in a compatible way and we transfer the assumptions (a)-(b) of Theorem \ref{mainpreciso} to the these stratifications;
  \item Using the assumption $(a)$ of Theorem \ref{mainpreciso} and the spectral sequence for a stratified space, we construct a chain complex $^\mathbb R E_1^{\bullet,0}$ whose cohomology computes the Betti numbers of $X_{0,1}^{log}(\mathbb R)$;
  \item Using the assumption $(b)$ of Theorem \ref{mainpreciso} and the spectral sequence for equivariant cohomology (combined with the assumption $(a)$), we construct a filtration $F^{\bullet}_{i}\subseteq$ $ ^\mathbb R E_1^{\bullet,0}$ such that $F^{\bullet}_{i}/F^{\bullet}_{i+1}$ is isomorphic to a complex depending only on the complex geometry of $X\rightarrow C$;
  \item We conclude the proof using the spectral sequence for a filtered complex.
    \end{enumerate}
      
   \subsection{Construction of \texorpdfstring{$X_{0,1}^{log}(\mathbb C)$}, and \texorpdfstring{$X_{0,1}^{log}(\mathbb R)$},}
   We endow $X\rightarrow C$ with the real logarithmic structure defined in Section \ref{subsec : semistablecase real} (see also Sections \ref{subsec: semistablecaselog} and \ref{subsec: semistablecasean}) and recall the commutative diagram
	\begin{center}
			 \begin{tikzcd}
		X_{0,1}^{log}(\mathbb C)\arrow{d}\arrow{r}\arrow[bend left]{rr}{\alpha}&	 X_0^{log}(\mathbb C)\arrow{d}\arrow{r}{\pi}& X_0(\mathbb C)\\
		\{1\}\arrow{r}&	 0^{\log}(\mathbb C)=S^1,
	\end{tikzcd}
		 \end{center} 
   compatible with the involutions and its real counterpart 
   	\begin{center}
			 \begin{tikzcd}
		X_{0,1}^{log}(\mathbb R)\arrow{d}\arrow{r}\arrow[bend left]{rr}{\alpha_{\mathbb R}}&	 X_0^{log}(\mathbb R)\arrow{d}\arrow{r}{\pi_{\mathbb R}}& X_0(\mathbb R)\\
		\{1\}\arrow{r}&	 0^{\log}(\mathbb R)=\{\pm 1\}.
	\end{tikzcd}
		 \end{center} 
By Theorems \ref{ogusnakayama} and \ref{katonakayamareal}, the general smooth fiber of $X(\mathbb C)\rightarrow C(\mathbb C)$ is homeomorphic to $X_{0,1}^{log}(\mathbb C)$ and the general positive smooth fiber near $0$ of $X(\mathbb R)\rightarrow C(\mathbb R)$ is homeomorphic to $X_{0,1}^{log}(\mathbb R)$. Hence, to prove Theorem \ref{mainpreciso}(1), we can replace $X_t(\mathbb R)$ with $X_{0,1}^{log}(\mathbb R)$. 
   \subsection{Stratifications}\label{subsec stratificationacaso}
   \numberwithin{equation}{subsection} 

   As in Sections \ref{subsec: stratification} and \ref{subsec : stratlogreal}, the collections 
   $$\{X^{0,log}_{\Delta,1}(\mathbb C):=\alpha^{-1}(X^0_{\Delta}(\mathbb C))\} \quad \text{and} \quad \{X^{0,log}_{\Delta,1}(\mathbb R):=\alpha_{\mathbb R}^{-1}(X^0_{\Delta}(\mathbb R))\} $$
   give equivariant stratifications of $X^{log}_{0,1}(\mathbb C)$ and $X^{log}_{0,1}(\mathbb R)$ respectively. Recall that, by Remarks \ref{keyproperty} and \ref{keypropertyreal}, one has that
   \begin{enumerate}
   	\item $X^{log}_{0,1}(\mathbb C)$ is a smooth $C^{\infty}$-manifold of dimension $2\Dim(X^0_{\Delta})+\vert \Delta\vert -1$;
   	\item $X^{log}_{0,1}(\mathbb R)$ is a smooth $C^{\infty}$-manifolds of dimension $\Dim(X^0_{\Delta})$.
   \end{enumerate}

   The following is the key lemma that allows us to transfer the assumptions (a) and (b) of Theorem \ref{mainpreciso} from the stratifications of $X_0$ to the ones of $X^{\log}_{0,1}(\mathbb C)$ and  $X^{\log}_{0,1}(\mathbb R)$.
  \begin{lemma}\label{lemma transfer assumptions}
  	\begin{enumerate}
  	\item[]
  	\item[(a)] $H^i(X^{0,\log}_{\Delta,1}(\mathbb R))=0$ for $i\geq 1$;
  	\item[(b)] $X^{0,\log}_{\Delta,1}(\mathbb C)$ is a maximal variety.
  	\end{enumerate}
  \end{lemma} 
   \proof
By Remark \ref{keypropertyreal}, one has that $X^{0,\log}_{\Delta,1}(\mathbb R)\rightarrow X^{0}_{\Delta}(\mathbb R)$ is a topological cover of degree $2^{\vert\Delta\vert-1}$. Since $H^1(X^0_{\Delta}(\mathbb R),\mathbb Z/2\mathbb Z)=0$ by assumption, this cover is the trivial one. Hence $X^0_{\Delta,1}(\mathbb R)$ is the disjoint union of $2^{\vert\Delta\vert-1}$ copies of $X^0_{\Delta}(\mathbb R)$, so that 
\begin{equation}\label{numeretti real}
b_i(X^{0,\log}_{\Delta,1}(\mathbb C))=0\quad  \text{for $i\geq 1$}\quad \text{and} \quad b_*(X^{0,\log}_{\Delta,1}(\mathbb C))=2^{\vert \Delta\vert -1}b_*(X^0_{\Delta}(\mathbb R)),
\end{equation}
and, in particular, we get (a). 

To prove (b), we first observe that, by the Smith-Thom inequality (\ref{eqn: smiththominequality}), one has 
$$b_*(X^{0,log}_{\Delta,1}(\mathbb R))\leq b_*(X^{0,log}_{\Delta,1}(\mathbb C)),$$
hence it is enough to prove that 
$$b_*(X^{0,log}_{\Delta,1}(\mathbb C))\leq b_*(X^{0,log}_{\Delta,1}(\mathbb R)).$$
First of all, one has the following relations
$$b_*(X^{0,log}_{\Delta,1}(\mathbb R))=2^{\vert \Delta\vert -1}b_*(X^0_{\Delta}(\mathbb R))=2^{\vert \Delta\vert -1}b_*(X^0_{\Delta}(\mathbb C)),$$
where the first equality follows from (\ref{numeretti real}) and the second one from the assumption (b) of Theorem \ref{mainpreciso} on $X^0_{\Delta}(\mathbb C).$ Hence, in order to conclude, it is enough show that 
$$b_*(X^{0,log}_{\Delta,1}(\mathbb C)) \leq 2^{\vert \Delta\vert -1}b_*(X^0_{\Delta}(\mathbb C)).$$
By Example \ref{exampleanilytification}(7) the morphism $\alpha:X^{0,\log}_{\Delta,1}(\mathbb C)\rightarrow X^0_{\Delta}(\mathbb C)$ is a locally trivial fibration with fiber $(S^1)^{\vert \Delta\vert -1}$. 
The Leray spectral sequence associated to such morphism $\alpha$
$$E_2^{a,b}:=H^a(X^{0}_{\Delta}(\mathbb C), R^b\alpha_*\Z/2\Z) \Rightarrow H^{a+b}(X^{0,\log}_{\Delta,1}(\mathbb C))$$
 shows that

$$b_i(X^{0,\log}_{\Delta,1}(\mathbb C))\leq \sum_a \Dim(H^a(X^{0}_{\Delta}(\mathbb C), R^{i-a}\alpha_*\Z/2\Z)).$$ 
By Lemma \ref{ilfascioecostante}, proved later, the locally constant sheaf $R^i\alpha_*\Z/2\Z$ is constant so that 
$$R^i\alpha_*\Z/2\Z\simeq (\Z/2\Z)^{\Dim(\displaystyle{H^i((S^1)^{\vert \Delta \vert -1})})}.$$
Hence 
$$\dim(H^a(X^{0}_{\Delta}(\mathbb C), R^{i-a}\alpha_*\Z/2\Z))=\dim(H^a(X^{0}_{\Delta}(\mathbb C), \Z/2\Z^{\dim(\displaystyle{H^{i-a}((S^1)^{\vert \Delta \vert -1})})}))=$$
$$=b_a(X^{0}_{\Delta}(\mathbb C))\dim(H^{i-a}((S^1)^{\vert \Delta \vert -1}))=b_a(X^{0}_{\Delta}(\mathbb C))\binom{\vert \Delta \vert -1}{i-a}.$$
This implies that
$$\sum_ib_i(X^{0,\log}_{\Delta,1}(\mathbb C))\leq \sum_i \sum_a b_a(X^{0}_{\Delta}(\mathbb C))\binom{\vert \Delta \vert -1}{i-a}=\sum_a \binom{\vert \Delta \vert-1}{a}\sum_i b_{i-a}(X^{0}_{\Delta}(\mathbb C))$$
$$\leq \sum_a \binom{\vert \Delta \vert -1}{a}\sum_i b_{i}(X^{0}_{\Delta}(\mathbb C))=  2^{\vert \Delta \vert  -1}\sum_a b_a(X^{0}_{\Delta}(\mathbb C)).$$
\endproof
\begin{lemma}\label{ilfascioecostante}
$R^i\alpha_*\Z/2\Z$ is a constant sheaf.
\end{lemma}
\proof
Let the notation be as in the following canonical commutative diagram
 \begin{center}
			 \begin{tikzcd}
			 X^{0,\log}_{\Delta,1}(\mathbb C)\arrow[hook]{r}{j}\arrow{rd}{\alpha} &  X^{0,\log}_{\Delta}(\mathbb C)\arrow[hook]{r}{h} \arrow{d}{\pi}  &  X_0^{\log}(\mathbb C)\arrow{d}\\
			 & X^{0}_{\Delta}(\mathbb C)\arrow[hook]{r} & X_0(\mathbb C).
			 	\end{tikzcd}
	 \end{center}
We first show that $R^i\pi_*\Z/2\Z$ is constant. By \cite[Lemma 1.5]{KatNak99}, one has that $R^i\pi_*\Z/2\Z\simeq (\bigwedge^i\overline M_{X^{0,\log}_{\Delta}}^{gp})\otimes \Z/2\Z$, so that it is enough to show that $\overline M_{X^{0,\log}_{\Delta}}^{gp}$ is constant. By \cite[Lemma 1.8.1]{Nak2000} one has 
$$\overline M_{X_0^{\log}}^{gp}\simeq \bigoplus_{i\in I} (X_i\rightarrow X)_*\mathbb Z$$
so that 
$$\overline M_{X^{0,\log}_{\Delta}}^{gp}=h^*\overline M_{X_0^{\log}}^{gp}=\bigoplus_{i\in I} h^*(X_i\rightarrow X)_*\mathbb Z= \bigoplus_{i\in I \text{}\vert \text{} X^0_{\Delta}\subseteq X_i}\mathbb Z.$$

We now deduce that this implies that $R^i\alpha_*\Z/2\Z$ is constant. Since $j$ is a closed immersion $R^i\alpha_*\Z/2\Z=R^i\pi_*(j_*\Z/2\Z)$, therefore it is enough to show that $R^i\pi_*\Z/2\Z\rightarrow R^i\pi_*(j_*\Z/2\Z)$ is surjective. This can be checked locally on $X^{0}_{\Delta}(\mathbb C)$. By Example \ref{exampleanilytification}(4), we can assume that $X^{0,log}_{\Delta}(\mathbb C)\simeq (S^1)^{\vert\Delta \vert}\times X^{0}_{\Delta}(\mathbb C)$ and $X^{0,log}_{\Delta,1}(\mathbb C)$ is the subset defined by the product of the first $\vert\Delta \vert$ components equal to $1$. In this case the result follows from the fact that the inclusion $X^{0,log}_{\Delta,1}(\mathbb C)\subseteq X^{0,log}_{\Delta}(\mathbb C)$ has a (non-canonical) retraction, depending on the choice of one of the first $\vert\Delta \vert$ components. 
\endproof
   \subsection{Computations of the real Betti numbers}\label{computationrealnumber}
 
  % Let $\omega: P \rightarrow \Z$ be the function that sends $\Delta$ to $-\Dim(X^0_\Delta)$. Since $\omega(\Delta)<\omega(\Delta')$ if $X^{0,log}_{\Delta,1}(\mathbb R)<X^{0,log}_{\Delta',1}(\mathbb R)$, namely $X^{log}_{\Delta,1}(\mathbb R)\supset X^{log}_{\Delta',1}(\mathbb R)$, 
   
  Let us consider the natural spectral sequence (see e.g. \cite{Pete16} and \cite[(3) pag. 2527]{Pete17})
   \begin{equation}\label{spectral sequence stratification real}
   	^{\R}E_1^{a,b}:=\bigoplus_{\Dim(X^0_\Delta)=a} H_c^{a+b}(X^{0,log}_{\Delta,1}(\mathbb R))\Rightarrow H_c^{a+b}(X_{0,1}^{log}(\mathbb R)).
   \end{equation}
Since $X^{0,log}_{\Delta,1}(\mathbb R)$ is a smooth $C^{\infty}$-manifolds of dimension $\Dim(X^0_{\Delta})$, we can apply Poincar\'e duality to deduce that
$$\bigoplus_{\Dim(X^0_\Delta)=a} H_c^{a+b}(X^{0,log}_{\Delta,1}(\mathbb R)) \simeq \bigoplus_{\Dim(X^0_\Delta)=a} H^{-b}(X^{0,log}_{\Delta,1}(\mathbb R)).$$
 By Lemma \ref{lemma transfer assumptions}(a), we get that $^{\R}E_1^{a,b}=0$ if $b\neq 0$. Hence the first page of the spectral sequence (\ref{spectral sequence stratification real}) reduces to the line $^\mathbb R E_1^{\bullet,0}$:
\begin{equation}\label{cochain complex}
\small{\dots \rightarrow\bigoplus_{\Dim(X^0_\Delta)=a-1} H_c^{a-1}(X^{0,log}_{\Delta,1}(\mathbb R))\xrightarrow{d_{a-1}} \bigoplus_{\Dim(X^0_\Delta)=a}H_c^{a}(X^{0,log}_{\Delta,1}(\mathbb R))\xrightarrow{d_{a}} \bigoplus_{\Dim(X^0_\Delta)=a+1} H_c^{a+1}(X^{0,log}_{\Delta,1}(\mathbb R))\rightarrow \dots}
\end{equation}
so that from (\ref{spectral sequence stratification real}) and Theorem \ref{katonakayamareal}, one gets
\begin{equation}\label{calcolotutto}
	H_c^a(X_t(\mathbb R))=H^a(\text{}^\mathbb R E_1^{\bullet,0}),
\end{equation}
for $t$ positive and near $0$. Hence, in order to conclude the proof of Theorem \ref{mainpreciso},(1) it is enough to bound the dimension of $H^a(\text{}^\mathbb R E_1^{\bullet,0})$.
   \subsection{Construction of the filtration}\label{subsec : constructionfiltration}
   Since $X^{0,log}_{\Delta,1}(\mathbb C)$ and $X^{0,log}_{\Delta,1}(\mathbb R)$ are smooth $C^{\infty}$-varieties,  by Poincar\'e duality and Lemma \ref{lemma transfer assumptions}(b), one has that $X^{0,log}_{\Delta,1}(\mathbb C)$ is maximal. The dimension of $X^{0,log}_{\Delta,1}(\mathbb R)$ is $\dim(X^0_{\Delta})$, therefore we have
   $$H_c^{\Dim(\Delta)}(X^{0,log}_{\Delta,1}(\mathbb R)))=H_c^*(X^{0,log}_{\Delta,1}(\mathbb R)).$$
   
   In order to get a filtration of $H_c^{\Dim(X^0_{\Delta})}(X^{0,log}_{\Delta,1}(\mathbb R))$, we use the tools of Section \ref{filtrazionest}. Once we set $m=2\dim(X_0)+1$ and choose $N\gg m$, one gets (using the construction of Section \ref{filtrazionest})  a decreasing filtration $\text{}_{m}F_{\Delta,\bullet}$ 
   $$\dots \subseteq  \text{}_{m}F_{\Delta,i+1} \subseteq\text{}_{m}F_{\Delta,i}\subseteq \dots  \subseteq \text{}_{m}F_{\Delta,0}=H^{\dim(X^0_\Delta)}_c(X^{0,log}_{\Delta,1}(\mathbb R))\big (\simeq H_{G,c}^{m}(X^{0,log}_{\Delta,1}(\mathbb C)\big )$$
  with graded quotients
$$\Gr_{_mF_{\Delta,\bullet}}^{i}:=\text{}_{m}F_{\Delta,i}/_{m}F_{\Delta,i+1}\simeq H_c^{m-i}(X^{0,log}_{\Delta,1}(\mathbb C))$$
Let $X^0_{\Delta}$ and $X^0_{\Delta'}$ be two strata such that $\Dim(X^0_{\Delta})=a-1=\Dim(X^0_{\Delta'})-1$ and $X^0_{\Delta}\subseteq X_{\Delta'}$ so that there is a commutative diagram
\begin{center}
			 \begin{tikzcd}
H_c^{m}(X^{0,\log}_{\Delta, 1}(\mathbb C)_{G})\arrow{r} & H_c^{m+1}(X^{0,\log}_{\Delta', 1}(\mathbb C)_{G})\\

H_c^{m}(X^{0,\log}_{\Delta, 1}(\mathbb R)\times BG)\arrow{r}\arrow{u}{\simeq} & H_c^{m+1}(X^{0,\log}_{\Delta', 1}(\mathbb R)\times BG)\arrow{u}{\simeq}\\

\underset{k}{\bigoplus} H_c^{m-k}(X^{0,\log}_{\Delta, 1}(\mathbb R))\otimes H_c^{k}(BG)\arrow{r}\arrow{u}{\simeq} & \underset{k}{\bigoplus} H_c^{m+1-k}(X^{0,\log}_{\Delta', 1}(\mathbb R))\otimes H_c^{k}(BG) \arrow{u}{\simeq}\\

H_c^{a-1}(X^{0,\log}_{\Delta, 1}(\mathbb R))\otimes H_c^{m-a+1}(BG)\arrow{u}{\simeq}\arrow{r} & H_c^{a}(X^{0,\log}_{\Delta', 1}(\mathbb R))\otimes H_c^{m-a+1}(BG)\arrow{u}{\simeq}\\
H_c^{a-1}(X^{0,\log}_{\Delta, 1}(\mathbb R))\arrow{u}{\simeq}\arrow{r}& H_c^{a}(X^{0,\log}_{\Delta', 1}(\mathbb R)),\arrow{u}{\simeq}
	\end{tikzcd}
	 \end{center}
	 where, for $\# \in\{\Delta, \Delta'\}$, 
	\begin{itemize}
		\item the horizontal maps 	are induced from the stratifications $\{X^{0,\log}_{\Delta, 1}(\mathbb C)_{G}\}$, $\{X^{0,\log}_{\Delta, 1}(\mathbb R)\times BG\}$ and $\{X^{0,\log}_{\Delta, 1}(\mathbb R)\}$ of $X^{\log}_{0,1}(\mathbb C)_{G}$,  $X^{\log}_{0,1}(\mathbb R)\times BG$ and  $X^{\log}_{0,1}(\mathbb R)$ respectively;
		\item the upper vertical maps are induced by the inclusions $X^{0,\log}_{1,\#}(\mathbb R)\times BG\subseteq X^{0,\log}_{1,\#}(\mathbb C)_{G}$ and they are isomorphisms by \cite[IV, Pag. 55, Application 3.7(b)]{borel} (see Section \ref{subsubsec maximalityequivariant});
		\item the remaining top-to-bottom vertical maps are induced respectively by the Kunneth formula, the fact that $H_c^{i}(X^{0,\log}_{1,\#}(\mathbb R))=0$ for $i\neq \Dim(X^0_{\Delta})$ and $H_c^{m-a+1}(BG)\simeq \Z/2\Z$.
	\end{itemize} 
Since the morphisms $H_c^{i}(X^{0,\log}_{\Delta, 1}(\mathbb C)_{G})\rightarrow H_c^{i+1}(X^{0,\log}_{\Delta', 1}(\mathbb C)_{G})$ and $H_c^{i}(X^{0,\log}_{\Delta, 1}(\mathbb C))\rightarrow H_c^{i+1}(X^{0,\log}_{\Delta', 1}(\mathbb C))$ give rise to a morphism of spectral sequences
\begin{center}
			 \begin{tikzcd}
  	^{ST}_{X^{0,\log}_{\Delta, 1}(\mathbb C)}E_2^{a,b}:=H^a(G,H_c^b(X^{0,\log}_{\Delta, 1}(\mathbb C)))\arrow[Rightarrow]{r}\arrow{d} & H_c^{a+b}(X^{0,\log}_{\Delta, 1}(\mathbb C)_{G})\arrow{d} \\
  		^{ST}_{X^{0,\log}_{\Delta', 1}(\mathbb C)}E_2^{a,b+1}:=H^a(G,H_c^{b+1}(X^{0,\log}_{\Delta', 1}(\mathbb C)))\arrow[Rightarrow]{r} & H_c^{a+b+1}(X^{0,\log}_{\Delta', 1}(\mathbb C)_{G}),
	\end{tikzcd}
	 \end{center}
the morphism $H^{a-1}_c(X^{0,log}_{\Delta,1}(\mathbb R))\rightarrow H^{a}_c(X^{0,log}_{\Delta',1}(\mathbb R))$ sends $_{m}F_{\Delta,i}$ to $_{m+1}F_{\Delta',i}$.

Therefore, if we set 
$$_mF^{p}_i:=\bigoplus_{\Dim(X^0_{\Delta})=p} \text{   }_{m+p}F_{\Delta,i}$$
we get a decreasing filtration $F^{\bullet}_{\bullet}$ of $^\mathbb R E_1^{\bullet,0}$ (Section \ref{computationrealnumber})
\begin{center}
			 \begin{tikzcd}[column sep=small]
			 _mF^{a-1}_{i+1}\arrow{r}\arrow[hook]{d} & _mF^{a}_{i+1}\arrow[hook]{d}  \arrow{r} &_mF^{a+1}_{i+1}\arrow[hook]{d} \\
	 _mF^{a-1}_{i}\arrow{r}\arrow[hook]{d} &_mF^{a}_{i}\arrow[hook]{d}  \arrow{r} &_mF^{a+1}_{i}\arrow[hook]{d} \\
			 			\dots\arrow{r}\arrow[hook]{d} &\dots\arrow[hook]{d}  \arrow{r} &\dots\arrow[hook]{d} \\
			 
	\dots \underset{\Dim(X^0_\Delta)=a-1}\bigoplus H_c^{a-1}(X^{0,log}_{\Delta,1}(\mathbb R))\arrow{r} & \underset{\Dim(X^0_\Delta)=a}\bigoplus H_c^{a}(X^{0,log}_{\Delta,1}(\mathbb R))\arrow{r}&\underset{\Dim(X^0_\Delta)=a+1}\bigoplus H_c^{a+1}(X^{0,log}_{\Delta,1}(\mathbb R)) \dots 
	\end{tikzcd}
	 \end{center}
	such that the graded quotient 
	$$\Gr_{F_{\bullet}}^{i, \bullet}:=\text{ }_mF^{\bullet}_{i}/\text{ }_mF^{\bullet}_{i+1}$$
	 is the complex $^{\C}E_1^{\bullet,m-i}$ of Section \ref{subsec: construction}
	 \begin{center}
			 \begin{tikzcd}[column sep=small]
			 _mF^{a-1}_{i}/_mF_{i+1}^{a-1}\arrow{r}\arrow{d}{\simeq} & _mF^{a}_{i}/_mF_{i+1}^{a}\arrow{r} \arrow{d}{\simeq}& _mF^{a+1}_{i}/_mF_{i+1}^{a+1}\arrow{d}{\simeq}\\
			 \underset{\Dim(X^0_\Delta)=a-1}\bigoplus H_c^{m+a-i-1}(X^{0,log}_{\Delta,1}(\mathbb C))\arrow{r} & \arrow{r}\underset{\Dim(X^0_\Delta)=a}\bigoplus H_c^{m+a-i}(X^{0,log}_{\Delta,1}(\mathbb C))\arrow{r}&
			 \underset{\Dim(X^0_\Delta)=a+1}\bigoplus H_c^{m+a-i+1}(X^{0,log}_{\Delta,1}(\mathbb C)).
			 	\end{tikzcd}
	 \end{center}	
   \subsection{End of the proof}
Recall the spectral sequence for a filtered complex
$$^FE_1^{r,s}:=H^{r+s}(\Gr_{F_{\bullet}}^{r, \bullet})\simeq H^{r+s}(\text{}^{\C}E_1^{\bullet,m-r})\Rightarrow H^{r+s}( ^\mathbb R E_1^{\bullet,0}).$$
Since 
$$H^{r+s}( \text{}^\mathbb R E_1^{\bullet,0})\simeq H^{r+s}(X_{0,1}^{\log}(\mathbb R)) \simeq H^{r+s}(X_t(\mathbb R)),$$
where the first isomorphism follows from (\ref{calcolotutto}) and second from Theorem \ref{katonakayamareal} (and the compactness of $X_t(\mathbb R)$), we get that
$$\Dim(H^{p}(X_t(\mathbb R),\Z/2\Z))=\sum_q \Dim(\text{}^FE_{\infty}^{q,p-q})\leq \sum_q \Dim(\text{}^FE_{1}^{q,p-q})=\sum_q \Dim(H^{p}(\text{}^{\C}E_1^{\bullet,m-q}))=\sum_q \Dim(H^{p}(\text{}^{\C}E_1^{\bullet,q})).$$

 \section{Proof of Theorem \texorpdfstring{\ref{mainpreciso}}-(2)}  \label{sec: hodge_number}
 	\numberwithin{equation}{subsection} 
 \subsection{Notation}
  Retain the notation and the assumption of Theorem \ref{mainpreciso}(2).
  We set 
  $$X^{\log}_{0,1}(\mathbb C):=X^{\log}_{0,1}\quad , \quad X^{0,\log}_{\Delta,1}(\mathbb C):=\widetilde {X}^{0}_{\Delta}$$
and for every topological space $Z$ we write $H^i(Z):=H^{i}(Z,\mathbb Q)$. 
Recall that there are morphisms
\begin{equation}\label{eq : morphism}
	\pi:X^{\log}_{0,1}\rightarrow X_0 \quad \text{and}\quad \pi_{\Delta}:\widetilde X^0_{\Delta}\rightarrow X^0_{\Delta}
\end{equation}
and stratifications
\begin{equation}\label{eq : stratifications}
\coprod X^0_{\Delta}=X_0 \quad \text{and}\quad \coprod \widetilde X^0_{\Delta}=X^{log}_{0,1}.
\end{equation}
For a $\Q$-mixed Hodge structure $V$, we write $\Gr^W_{i}$ for its weight $i$-graded piece.

Fix $q\in \mathbb N$.
  \subsection{Strategy}\label{sec : strategy}
  We detail the strategy to prove Theorem \ref{mainpreciso}(2).
  \begin{enumerate}
  	\item In order to prove point $(ii)$, what one needs to do is to transfer the assumption on the torsion freeness on the cohomology of $X^0_{\Delta}$ to the one of $\widetilde X^0_{\Delta}$. This reduces to study the Leray spectral sequence with integer coefficient for the morphism
  $\pi_{\Delta}: \widetilde X^0_{\Delta}\rightarrow X^0_{\Delta}$ and to show that it degenerates in $E_2$ (Lemma \ref{cor: torsion_free}).
  \item We then deal with point $(i)$. This involves a careful analysis of different spectral sequences coming from geometry. Following \cite{tropicalhomology}, we first show that the $(p,q)$-hodge number of the general fiber can be computed as the weight $2q$-part of the limiting Hodge structure on $H^{p+q}(X^{\log}_{0,1})$ (Lemma \ref{lem : hodgenumberweight});
  \item The Leray spectral sequence $^LE_2$ for $X^{\log}_{0,1}\displaystyle{\stackrel{\pi}{\longrightarrow}} X_0$ is a spectral sequence of mixed Hodge structure which degenerates in $^LE_3$, hence $\Gr^W_{2q}H^{p+q}$ can be computed via its second page (Section \ref{sec : vanishingspectral});
  \item We then filter the complex $\text{}^{\mathbb C}E_1^{\bullet, q}$ by using the filtration induced on every term by the Leray spectral sequence for $\widetilde X^0_{\Delta}\displaystyle{\stackrel{\pi_{\Delta}}{\longrightarrow}}  X^0_{\Delta}$ and we consider the spectral sequence $^FE_1$ for the filtered complex (Section \ref{sec : filtration via leray});
  \item Exploiting the assumptions on weights, we show that the first page of $^FE_1$ is isomorphic, up to change the indexing, to the second page of $^LE_2$. Hence the cohomology of the lines of $^FE_1$ computes the weights of the limiting Hodge structure (Section \ref{sec : gigacomputation});
  \item This arguments bounds the cohomology of $\text{}^{\mathbb C}E_1^{\bullet, q}$ with the Hodge numbers of the generic fiber and to get an actual equality, one has to do a few computations recalling how $\text{}^{\mathbb C}E_1^{\bullet, q}$ was constructed (Section \ref{sec : preliminary reductions}). 
  	  \end{enumerate}
\begin{remark}
Another natural strategy to prove Theorem \ref{mainpreciso}(ii) is to try to reduce to finite fields, by the specialization arguments from \cite{Nak2000}, and to use the theory of Frobenius weights there. While this might be possible, one issue is that the existence of the spectral sequence with \'etale cohomology with compact support for logarithmic schemes has not been proved. While this spectral sequence seems likely to exist and it seems to be useful to study it, doing it would take us far away from the techniques used in the rest of this paper. Hence we preferred to remain in the complex analytic setting, at the cost of doing a somehow more involved argument. We hope to come back to the \'etale case problem in the near future.
\end{remark}
Before starting with the actual proof of Theorem \ref{mainpreciso}(ii) in Sections \ref{sec : preliminary reductions} and \ref{sec : proof of (2)}, we collect some preliminaries on the nearby cycles functors in Section \ref{nearby cycles}.
\subsection{Nearby cycles and spectral sequences}\label{nearby cycles}
\numberwithin{equation}{subsubsection} 

\subsubsection{Computation of nearby cycles over open strata}
For the lack of a reference, we collect a few basic results on the nearby cycle functors $R^i\pi_{\Delta,*}A$ for an abelian group $A$. This is essentially a more refined version of Lemma \ref{ilfascioecostante}. The analogue description in the \'etale case has been done in \cite[Appendix A.1]{Nak2000}. 

First, observe that $\pi_{\Delta,*}A\simeq A$, so that we want to compute $R^i\pi_{\Delta,*}A$ for $i\geq 1$. Recall that $\pi_{\Delta}:\widetilde X_{\Delta,1}\rightarrow X^{0}_{\Delta}$ factorises as 
$$\widetilde X^0_{\Delta}\xrightarrow{\alpha_{\Delta}} X^{0,log}_{\Delta} \xrightarrow{\beta_{\Delta}}X^{0}_{\Delta}.$$
Hence, since $\alpha_{\Delta}$ is a closed immersion,
$$R^1\pi_{\Delta*}A\simeq R^1\beta_{\Delta,*}\alpha_{\Delta,*}A.$$
Moreover, since $\pi_{\Delta}:\widetilde X^{0}_{\Delta}\rightarrow X^{0}_{\Delta}$ is a locally trivial $(S^1)^{\vert \Delta \vert -1}$-bundle, the natural cup product map
$$\wedge^i R^1\pi_{\Delta,*}A\xrightarrow{\simeq} R^i\pi_{\Delta,*}A$$
is an isomorphism.

Similar reasonings show that 
\begin{equation}\label{eq : wedge property}
	\wedge^i R^1\beta_{\Delta,*}A\xrightarrow{\simeq} R^i\beta_{\Delta,*}A
\end{equation}
is an isomorphism and, as in Lemma \ref{ilfascioecostante}, by \cite[Lemma 1.5]{KatNak99} and \cite{Nak2000},
there is a canonical isomorphism
\begin{equation}\label{eq : computationR1}
	R^1\beta_{\Delta,*}A\simeq \bigoplus_{X^0_\Delta\subseteq X_i} A(-1).
\end{equation}
Working locally as in Lemma \ref{ilfascioecostante}, one sees that the natural morphism
$$R^1\beta_{\Delta,*}A\rightarrow R^1\pi_{\Delta,*}A\simeq R^1\beta_{\Delta,*}\alpha_{\Delta,*}A$$
is surjective and fits into an exact sequence
$$0\rightarrow A(-1)\xrightarrow{diag}  R^{1}\beta_{\Delta,*}A\simeq\bigoplus_{X^0_\Delta\subseteq X_i} A(-1) \rightarrow  R^{1}\pi_{\Delta,*}A\rightarrow 0$$
where $diag$ is the diagonal map. Observe that the choice of an $i$ such that $X^0_\Delta\subseteq X_i$, induces a (non-canonical) splitting of the exact sequences.
Putting all together, we get a natural exact sequence
$$0\rightarrow A(-1)\otimes R^{i-1}\pi_{\Delta,*}A\rightarrow  R^{i}\beta_{\Delta,*}A\rightarrow  R^{i}\pi_{\Delta,*}A\rightarrow 0$$
which splits (non-canonically).
In particular there is a canonical short exact sequence
\begin{equation}\label{eq : splitexact}
	0\rightarrow H_c^j(X^0_{\Delta},A(-1)\otimes R^{i-1}\pi_{\Delta,*}A)\rightarrow  H_c^j(X^0_{\Delta},R^{i}\beta_{\Delta,*}A)\rightarrow  H_c^j(X^0_{\Delta},R^{i}\pi_{\Delta,*}A)\rightarrow 0,
\end{equation}
which splits (non canonically).
\subsubsection{Weights and Hodge structures on open strata}
Now by (\ref{eq : computationR1}) and (\ref{eq : wedge property})
$$H_c^j(X^0_{\Delta},R^{i}\beta_{\Delta,*}A)\simeq H_c^i(X^0_{\Delta},\wedge^iR^{1}\beta_{\Delta,*}A)\simeq H_c^i(X^0_{\Delta},\wedge^i(\bigoplus_{X^0_\Delta\subseteq X_i} A(-1)))\simeq H_c^j(X^0_{\Delta},A(-i))^{\binom{\vert \Delta \vert -1}{i}}.$$
Combining this with (\ref{eq : splitexact}) and Poincar\'e Duality we get the following.
\begin{lemma}\label{lem : weights}
If $X^0_{\Delta}$ has dimension $d$, then  
$H_c^j(X^0_{\Delta},R^{i}\pi_*\Q)$ carries a natural Hodge structure pure of weight $2j-2d+2i$
\end{lemma}
We now apply Lemma \ref{lem : weights} to the Leray spectral sequences (with compact support) 
\begin{equation}\label{eq : spectral sequences leray open}
	\quad _{\Delta}^LE^{a,b}_2:=H_c^a(X^{0}_{\Delta},R^b\pi_{\Delta,*}\Q)\Rightarrow H_c^{a+b}(\widetilde X^{0}_{\Delta}).
\end{equation}
associated to the morphisms $\pi_{\Delta}: \widetilde X^0_{\Delta} \rightarrow X^0_{\Delta}$.
\begin{lemma}\label{lem : localLeraydegenerates in E2}
	$_{\Delta}^LE^{a,b}_2$ degenerates in $E_2$.
\end{lemma}
\proof
By Lemma \ref{lem : weights}, the term $_{\Delta}^LE^{a,b}_2$  is pure of weight $2a-2\Dim (X_{\Delta}^0)+2b$.  Since there are no morphism between mixed Hodge structures of different weights and the morphisms are compatible with the Hodge structures, the maps $_{\Delta}^LE^{a,b}_r \rightarrow _{\Delta}^LE^{a+r,b-r+1}_r$ are trivial for all $a,b,r\in \mathbb Z$.
\endproof
\subsubsection{Vanishing cycles spectral sequence}\label{sec : vanishingspectral}
Consider now the Leray spectral sequence for the morphism $\pi:X_{0,1}^{log}\rightarrow X_0$ 
\begin{equation}\label{eq : Leray spectral chiusa}
	^LE^{a,b}_2:=H^a(X_0,R^b\pi_*\Q)\Rightarrow H^{a+b}(X^{\log}_{0,1}).
\end{equation}
By \cite[Corollary 8.4]{illusielog}, the sheaf $R^b\pi_*\Q$ identifies with the classical sheaf of vanishing cycles $R^b\Psi \Q$. Hence, by \cite[Section 1.4]{DS04}, it is a spectral sequence of mixed Hodge structures, called also the the vanishing cycles spectral sequence.
By \cite{Sc1973} (see also \cite{SZ90} and the discussion in \cite[Section 8.8, before (8.8.6)]{illusielog}), it degenerates in $E_3$. 

Since taking graded pieces for the weight filtration is an exact operation, we can consider the weight $2q$-graded piece of the spectral sequence (\ref{eq : Leray spectral chiusa})
\begin{equation}\label{eqn:weight_2q_leray}
	\Gr^W_{2q}\text{ }^LE^{a,b}_2:=\Gr^W_{2q}H^a(X_0,R^b\pi_*\Q)\Rightarrow \Gr^W_{2q}H^{a+b}(X^{\log}_{0,1})
\end{equation}
which still degenerates in $E_3$, so that 
\begin{equation}\label{eq : weight in spectral}
	\Gr^W_{2q}H^{p+q}(X^{\log}_{0,1})=\sum_{i} \Gr^W_{2q}\text{ }^LE^{i,p+q-i}_3.
\end{equation}

\subsection{Preliminary reductions}\label{sec : preliminary reductions}
\subsubsection{Using torsion freeness}
We first prove point $(ii)$ in Theorem \ref{mainpreciso}(2). This is a consequence of the following corollary of Lemma \ref{lem : localLeraydegenerates in E2}, which transfers part of the $(c)$-assumptions from $X^{0}_{\Delta}$ to $\widetilde X^{0}_{\Delta}$.
\begin{corollary}
\label{cor: torsion_free}
$H_c^{n}(\widetilde X^{0}_{\Delta},\mathbb Z)$ is torsion free.
\end{corollary}
\proof
Let us consider the rational and the integer Leray spectral sequences $_{\mathbb Q}^LE_2$ and $_{\mathbb Z}^LE_2$ associated to the morphism $\pi_{\Delta}: \widetilde X^{0}_{\Delta} \rightarrow  X^{0}_{\Delta}$. Thanks to Lemma \ref{lem : localLeraydegenerates in E2}, we know that $_{\mathbb Q}^LE_2$ degenerates in $E_2$. On the other hand, by hypothesis of Theorem \ref{mainpreciso} and the split exact sequence (\ref{eq : splitexact}), the abelian group $_{\mathbb Z}^LE_2^{i,n-i}$ has no torsion for all $i$. Therefore the differential $_{\mathbb Z}^LE_r^{a,b} \rightarrow _{\mathbb Z}^LE_r^{a+r,b-r+1}$ has to be trivial for all $a,b,r\in \mathbb N$, so that $^L_{\mathbb Z}E_2$ degenerates in $E_2$. Hence, there exists a decreasing filtration $F^s \subseteq F^{s-1} \subseteq \dots \subseteq H_c^{n}(\widetilde X^{0}_{\Delta},\mathbb Z)$ such that $F_i/F_{i+1} \simeq \text{ }_{\mathbb Z}^LE_2^{i,n-i}$ has no torsion. This implies that $H_c^{n}(\widetilde X^{0}_{\Delta},\mathbb Z)$ has no torsion and the proof of the lemma is concluded. \endproof

Thanks to Corollary \ref{cor: torsion_free}         
$^{\mathbb C}_{\mathbb Z}E^{\bullet, q}_1\otimes \mathbb Z/2\mathbb Z\simeq \text{ }_{\mathbb Z/2\mathbb Z}^{\quad \mathbb C}E^{\bullet, q}_1$ and point $(ii)$ of Theorem \ref{mainpreciso}(2) is proved.
\subsubsection{Reduction to a weight computation}
We now start the proof of (i), by reducing it to a proof of the inequality (\ref{sottoassunzione}) in the following Lemma.
\begin{lemma}\label{lem : reduction}
	Assume that for all $p,q\in\mathbb N$ the inequality
	\begin{equation}\label{sottoassunzione}
	\Dim(H^p(\text{ }^{\mathbb C}E^{\bullet,q}_1))\leq \Gr^W_{2q}H^{p+q}(X_{0,1}^{log}).
\end{equation} 
holds. Then one has
		$\Gr^W_{2q}H^{p+q}(X_{0,1}^{log})=h^{p,q}$. 
\end{lemma}
In order to prove Lemma \ref{lem : reduction}, we start relating the Hodge numbers of the general fiber to the weight filtration in the cohomology of $X^{\log}_{0,1}$.

\begin{lemma}\label{lem : hodgenumberweight}
One has the equality $h^{p,q}(X_t)=\dim(\Gr^W_{2q}H^{p+q}(X^{\log}_{0,1}))$ and the mixed Hodge structure on $H^{p+q}(X^{\log}_{0,1})$ has only even weights.
	\end{lemma}
\proof
This is essentially proved in \cite[Middle of Page 31]{tropicalhomology}. We briefly recall the argument. For every subset $J\subseteq I$ consider the spectral sequence
$$^{\mathbb C} E^{a,b}_1:=\bigoplus_{\displaystyle{\substack{X^0_{\Delta}\subseteq X_J,\\ \Dim(X^0_{\Delta})=a}}} H_c^{a+b}(X^{0}_{\Delta})\Rightarrow H^{a+b}(X_J).$$
By assumption and Poincar\'e duality, $^\mathbb CE^{a,b}_1$ 
 is pure of type $(b,b)$, 
 hence $H^n(X_J)=0$ for $n$ odd and $H^{2n}(X_J)\simeq \mathbb Q(-n)^{n_J}$ for some $n_J\in \mathbb Z$.
Hence, the weight spectral sequence
$$E_2^{i,j}:=\bigoplus_{a\geq \max{\{0,i\}}}\bigoplus_{\vert J\vert =2a-i}H^{2i+j-2a}(X_J)(i-a)\Rightarrow H^{i+j}(X^{\log}_{0,1})$$
shows that $W^{2i}=F_{p+q-i}$, where $W^{j}$ is the increasing weight filtration on $H^{p+q}(X^{\log}_{0,1})$ and $F_{j}$ the decreasing Hodge filtration. Hence
$$h^{p,q}(X_t)=h^{q,p}(X_t)=\Dim(F_p/F_{p+1})=\dim(W^{2q}/W^{2q-2})= \Dim(\Gr^W_{2q}H^{p+q}(X^{\log}_{0,1})).$$
\endproof
\proof[Proof of Lemma \ref{lem : reduction}]
Lemma \ref{lem : hodgenumberweight} implies that
$$\dim(H^{p+q}(X_{0,1}^{log}))=\sum_{i}\Gr^W_{2i}H^{p+q}(X_{0,1}^{log}).$$
Therefore, by the spectral sequence (\ref{eq : spectral complex}), one gets the inequality 
$$\sum_i\Dim(H^{i}(\text{ }^{\mathbb C}E^{\bullet,p+q-i}_1))\geq\sum_{i}\Gr^W_{2i}H^{p+q}(X_{0,1}^{log}).$$
The assumption (\ref{sottoassunzione}) implies that 
$$\Dim(H^p(\text{ }^{\mathbb C}E^{\bullet,q}_1))=\Gr^W_{2q}H^{p+q}(X_{0,1}^{log}).$$
Hence, again by Lemma \ref{lem : hodgenumberweight}, one has $\Gr^W_{2q}H^{p+q}(X_{0,1}^{log})=h^{p,q}(X_t)$.
\endproof
\subsection{End of the proof}\label{sec : proof of (2)}
By Lemma \ref{lem : reduction}, in order to end the proof of Theorem \ref{mainpreciso}, it remains to prove the inequality (\ref{sottoassunzione}).
\subsubsection{Filtering the complex}\label{sec : filtration via leray}
Recall the spectral sequences 
$$^{r}_qE^{a,b}_1:=\bigoplus_{\Dim(X^0_{\Delta})=a}H_c^{a+b}(X^0_{\Delta},R^{q-r}\pi_{\Delta,*}\Q)\Rightarrow H^{a+b}(X_{0},R^{q-r}\pi_*\Q)$$
associated to the stratifications (\ref{eq : stratifications}).
By Lemma \ref{lem : localLeraydegenerates in E2}, there exists a decreasing filtration $_nF_{i,\Delta}$ of $H_c^{n}(\widetilde X_{\Delta}^0)$ such that 
$_nF_{i,\Delta}/_nF_{i+1,\Delta}\simeq H_c^{i}(X_{\Delta}^0,R^{n-i}\pi_{\Delta,*}\mathbb Q)$. 

If $\Dim(X^0_{\Delta})$ and $\Dim(X^0_{\Delta'})$ are respectively $m-1$ and $m$, and $X^0_{\Delta}\subseteq X_{\Delta'}$, there is a morphism of spectral sequences
\begin{center}
			 \begin{tikzcd}
			 ^L_{\Delta}E_2^{a,b}:= H_c^a(X^{0}_{\Delta},R^b\pi_{\Delta,*}\Q) \arrow[Rightarrow]{r}\arrow{d} & H_c^{a+b}(\widetilde X^{0}_{\Delta})\arrow{d}\\
			 			 ^L_{\Delta'}E_2^{a+1,b}:= H_c^{a+1}(X^{0}_{\Delta'},R^b\pi_{\Delta',*}\Q) \arrow[Rightarrow]{r}& H_c^{a+b+1}(\widetilde X^{0}_{\Delta'}),
	\end{tikzcd}
	 \end{center}
so that the morphism $H^{n-1}_c(\widetilde X^{0}_{\Delta})\rightarrow H^{n}_c(\widetilde X^{0}_{\Delta'})$ sends $_{n}F_{\Delta,i}$ to $_{n+1}F_{\Delta',i+1}$.
Hence, setting 
$$^qF^p_i:=\bigoplus_{\dim(\widetilde{X}^0_{\Delta})=p}\text{ }_{p+q}F_{\Delta',i+p}\subseteq \bigoplus_{\dim(\widetilde{X}^0_{\Delta})=p}H_c^{p+q}(\widetilde{X}^0_{\Delta})=\text{ }^{\mathbb C}E_1^{p,q},  $$ 
we get a decreasing filtration 
$$^qF^{\bullet}_{\bullet} :=   \quad ^qF^{\bullet}_{i+1}\subseteq ^qF^{\bullet}_{i}\subseteq \dots \subseteq \text{ }^qF^{\bullet}_{0}= \text{ }^\mathbb CE^{\bullet,q}_1$$
such that 
$$\Gr_{\text{ }^qF_{\bullet}}^{i, \bullet}:=\text{ }^qF^{\bullet}_{i}/^qF^{\bullet}_{i+1}\simeq \text{}^{i}_qE_1^{\bullet,i}.$$
Hence the spectral sequence for a filtered complex reads
\begin{equation}\label{eq : Efiltration}
\text{}_q^FE_1^{a,b}=H^{a+b}(\Gr_{\text{ }^qF_{\bullet}}^{a, \bullet})=H^{a+b}(\text{}^{a}_qE_1^{\bullet,a})=\text{ }^{a}_qE^{a+b,a}_2\Rightarrow H^{a+b}(\text{}^{\C}E^{\bullet,q}_1).
\end{equation} 
\subsubsection{Computation of the spectral sequence of the filtered complex}\label{sec : gigacomputation}
By (\ref{eq : Efiltration})
\begin{equation}\label{eq : infinity}
	\dim(H^{p}(\text{}^{\mathbb C}E^{\bullet, q}))\leq \sum_i \dim(\text{ }_q^FE_{2}^{i,p-i}) .
\end{equation}
We now fully employ the assumption on the weights by using Lemma \ref{lem : weights}.
\begin{lemma}
\label{cor : annulazione}
For every $r\in \mathbb N$ the spectral sequence 
$$^{r}_qE^{a,b}_1:=\bigoplus_{\Dim(X^0_{\Delta})=a}H_c^{a+b}(X^0_{\Delta},R^{q-r}\pi_{\Delta,*}\Q)\Rightarrow H^{a+b}(X_{0},R^{q-r}\pi_*\Q),$$
degenerates in $E_2$ and it induces a natural isomorphism
$$^r_qE_2^{a,b}\simeq \Gr^W_{2b+2q-2r}H^{a+b}(X_{0},R^{q-r}\pi_*\Q).$$
\end{lemma}
\proof
Since $H_c^{a+b}(X^0_{\Delta},R^{q-r}\pi_{\Delta,*}\Q)$ is pure of weight $2b+2q-2r$ by Lemma \ref{lem : weights} and there are no morphisms between Hodge structures of different weights, the spectral sequence degenerates in $E_2$ and it induces the displayed isomorphism.
\endproof
By (\ref{eq : Efiltration}) one has 
\begin{equation}\label{eq : computationE1}
	\text{}^F_qE_1^{a,b}=\text{}^{a}_qE^{a+b,a}_2,
\end{equation}
which by Corollary \ref{cor : annulazione}, is  canonically isomorphic to 
$$\text{}^{a}_qE^{a+b,a}_2\simeq \Gr^W_{2q}H^{2a+b}(X_{0},R^{q-a}\pi_*\Q),$$
which by definition (see (\ref{eqn:weight_2q_leray})) is $\Gr^W_{2q}\text{ }^LE^{2a+b,q-a}_{2}$.
Hence the bijection $\varphi:\mathbb Z^{2}\rightarrow \mathbb Z^{2}$ sending $(x,y)$ to $(2x+y,q-x)$ (with inverse $(x,y)\mapsto  (q-y, 2y+x-2q)$), induces a natural isomorphism 
\begin{equation}\label{eq : superiso}
	^F_qE^{a,b}_1\simeq \Gr^W_{2q}\text{ }^LE^{2a+b,q-a}_2.
\end{equation}
So one gets
$$H^p(\text{ }^{\mathbb C}E_1^{\bullet,q})\leq \sum_i\dim(\text{ }_q^FE_2^{i,p-i})=\sum_i \Dim(\Gr^W_{2q}\text{ }^LE^{p+i,q-i}_3)=\sum_i \Dim(\Gr^W_{2q}\text{ }^LE^{p+i,q-i}_{\infty})=\Gr^W_{2q}H^{p+q}(X^{log}_{0,1}),$$
where, to summarize, the first inequality comes from (\ref{eq : Efiltration}), the second equality from (\ref{eq : superiso}),  the third one from the degeneration of $\Gr^W_{2q}\text{ }^LE$ in $E_3$ (Section \ref{sec : vanishingspectral}) and the last one from (\ref{eq : Leray spectral chiusa}).
Hence we proved the inequality (\ref{sottoassunzione}) and the proof of Theorem \ref{mainpreciso}(2) is concluded.
\bibliographystyle{alpha}
\bibliography{bibliografia.bib}
Emiliano Ambrosi, \textsc{Universit\'{e} de Strasbourg}
\par\nopagebreak
 \textit{eambrosi'at'unistra.fr}: \texttt{ Institut de Recherche Math\'{e}matique Avanc\'{e}e (IRMA)}\\
Matilde Manzaroli, \textsc{Universit\"{a}t T\"{u}bingen}\par\nopagebreak
 \textit{matilde.manzaroli'at'uni-tuebingen.de}: \texttt{Mathematisches Institut}\end{document}